\documentclass[11pt]{amsart}
\usepackage[margin=1in]{geometry}
\usepackage[toc,page]{appendix}
\usepackage{amsmath}
\usepackage{amscd}
\usepackage{amssymb}
\usepackage{amsfonts}
\usepackage{amsxtra}
\usepackage{amsthm}
\usepackage{enumitem}
\usepackage{tikz-cd}
\usepackage{mathrsfs}
\usepackage{hyperref}
\usepackage{chngcntr}
\usepackage{apptools}
\usepackage{indentfirst}
\usepackage{color}
\usepackage{multirow}
\usepackage{bm}

\AtAppendix{\counterwithin{lemma}{section}}
\AtAppendix{\counterwithin{mydef}{section}}
\AtAppendix{\counterwithin{prop}{section}}
\AtAppendix{\counterwithin{theorem}{section}}

\theoremstyle{plain}
\newtheorem*{theorem*}{Theorem}
\newtheorem{theorem}{Theorem}[section]
\newtheorem{prop}[theorem]{Proposition}
\newtheorem{lemma}[theorem]{Lemma}
\newtheorem*{lemma*}{Lemma}
\newtheorem{corollary}[theorem]{Corollary}
\newtheorem{conj}[theorem]{Conjecture}
\theoremstyle{definition}
\newtheorem{exmp}[theorem]{Example}
\newtheorem{mydef}[theorem]{Definition}
\newtheorem{rmk}[theorem]{Remark}
\newtheorem{construction}[theorem]{Construction}

\DeclareMathOperator{\NN}{\mathbb{N}}
\DeclareMathOperator{\ZZ}{\mathbb{Z}}
\DeclareMathOperator{\QQ}{\mathbb{Q}}
\DeclareMathOperator{\RR}{\mathbb{R}}
\DeclareMathOperator{\CC}{\mathbb{C}}
\DeclareMathOperator{\PP}{\mathbb{P}}
\DeclareMathOperator{\FF}{\mathbb{F}}

\DeclareMathOperator{\cE}{\mathcal{E}}
\DeclareMathOperator{\cF}{\mathcal{F}}
\DeclareMathOperator{\cG}{\mathcal{G}}
\DeclareMathOperator{\cO}{\mathcal{O}}

\DeclareMathOperator{\cK}{\mathcal{K}}

\DeclareMathOperator{\cL}{\mathcal{L}}
\DeclareMathOperator{\cP}{\mathcal{P}}
\DeclareMathOperator{\cM}{\mathcal{M}}
\DeclareMathOperator{\cD}{\mathcal{D}}
\DeclareMathOperator{\cV}{\mathcal{V}}
\DeclareMathOperator{\cR}{\mathcal{R}}
\DeclareMathOperator{\cS}{\mathcal{S}}
\DeclareMathOperator{\cW}{\mathcal{W}}

\DeclareMathOperator{\bv}{\textup{\textbf{v}}}
\DeclareMathOperator{\bw}{\textup{\textbf{w}}}

\DeclareMathAlphabet\mathbfcal{OMS}{cmsy}{b}{n}

\DeclareMathOperator{\Hom}{Hom}

\DeclareMathOperator{\NS}{NS}

\DeclareMathOperator{\Ext}{Ext}
\DeclareMathOperator{\ext}{ext}
\DeclareMathOperator{\Pic}{Pic}
\DeclareMathOperator{\ch}{ch}

\DeclareMathOperator{\gr}{gr}
\DeclareMathOperator{\Nef}{Nef}

\DeclareMathOperator{\Quot}{\text{Quot}}

\DeclareMathOperator{\Eff}{\text{Eff}}

\title[Brill-Noether for del Pezzo surfaces]{Brill-Noether and existence of semistable sheaves for del Pezzo surfaces}

\date{\today}
\begin{document}

\author[D. Levine]{Daniel Levine}
\address{Department of Mathematics, The Pennsylvania State University, University Park, PA 16802}
\email{dul190@psu.edu}

\author[S. Zhang]{Shizhuo Zhang}
\address{Max-Planck institute for mathematics, Vivatsgasse 7, 53111 Bonn, Germany}
\email{shizhuozhang@mpim-bonn.mpg.de}

\thanks{During the preparation of this article the first author was partially supported by the NSF FRG grant DMS-1664303}

\begin{abstract}
Let $X_m$ be a del Pezzo surface of degree $9-m$. When $m \leq 5$, we compute the cohomology of a general sheaf in $M(\bv)$, the moduli space of Gieseker semistable sheaves with Chern character $\bv$. We also classify the Chern characters for which the general sheaf in $M(\bv)$ is non-special, i.e. has at most one nonzero cohomology group. Our results hold for arbitrary polarizations, slope semistability, and semi-exceptional moduli spaces. When $m \leq 6$, we further show our construction of certain vector bundles implies the existence of stable and semistable sheaves with respect to the anti-canonical polarization.
\end{abstract}

\maketitle

\tableofcontents

\section{Introduction}
Let $(X,H)$ be a smooth del Pezzo surface with a polarization $H$,  and let $M_{X,H}(\bv)$ be an irreducible moduli space of Gieseker $H$-semistable sheaves of Chern character $\bv$. We address two fundamental problems in the theory of moduli of sheaves: (1) The calculation of the cohomology of a general sheaf in $M_{X,H}(\bv)$, and (2) the classification of Chern characters $\bv$ for which $M_{X,H}(\bv)$ is nonempty.

\subsection{Higher rank Brill-Noether theory} The higher rank analogue of Brill-Noether theory consists of finding a classification for non-special Chern characters (Definition \ref{non-special definition}) on a polarized variety $(X,H)$ and computing the codimension of the locus of special sheaves. The first part of higher rank Brill-Noether theory is already interesting and is the focus of the majority of this paper. For example, the classification of non-special Chern characters is necessary in the study of Le Potier's strange duality conjecture. Other applications have been found in the classification of globally generated Chern characters (\cite{CoskunHuizengaBrillNoetherGloballyGenerated}) and the classification of Chern characters whose moduli spaces $M_{X,H}(\bv)$ are nonempty (\cite{Hirzebruchexist}).

In this paper, we classify the non-special Chern characters for del Pezzo surfaces of degree at least $4$. Moduli spaces of rank one sheaves are well-understood, but behave differently from moduli spaces of sheaves of higher rank (see \cite{CoskunHuizengaWeakBrillNoether}). Therefore, we only consider moduli spaces of sheaves of rank at least $2$.

Let $X_m$ be a del Pezzo surface of degree $9-m$ obtained by blowing up $\PP^2$ at $m$ general points. Recall that there is a Weyl group $W_m$ that acts on $\Pic(X_m)$ (\cite{ManinWeylGroup}). Let $L \in \Pic(X_m)$ be the total transform of a line from $\PP^2$ and let $\cP_{L^W}(\bv)$ be the stack of sheaves prioritary with respect to every Weyl group translate of $L$, i.e. $\Ext^2(\cE,\cE(-\sigma(L))) = 0$ for all $\sigma \in W_m$ (Section \ref{Brill-Noether}). The stack $\cP_{L^W}(\bv)$ is irreducible by a theorem of Walter (\cite{WalterIrreducibility}). When the stack of slope semistable sheaves $\cM_{X,H}^{\mu ss}(\bv)$ is nonempty, the stack $\cP_{L^W}(\bv)$ contains $\cM_{X,H}^{\mu ss}(\bv)$ as an open dense substack (Section \ref{prioritary section}).

\begin{mydef}
\label{non-special definition}
Let $\bv \in K(X_m)$ be a Chern character on $X_m$. A coherent sheaf $\cE$ is \textit{non-special} if it has at most one nonzero cohomology group. We say $\bv$ is \textit{non-special} if there exists a non-special sheaf in $\cP_{L^W}(\bv)$.
\end{mydef}

The classification of non-special Chern characters for $\PP^2$ was worked out in \cite{GottscheHirschowitzBrillNoetherP2} and for Hirzebruch surfaces, including $\PP^1 \times \PP^1$ and $\FF_1$, in \cite{CoskunHuizengaBrillNoetherGloballyGenerated}. For del Pezzo surfaces and arbitrary blowups, partial results were obtained in \cite{CoskunHuizengaWeakBrillNoether} under the condition that the Chern character $\bv$ satisfies $\chi(\bv) = 0$. Our main result gives a classification of all non-special Chern characters.

\begin{theorem}[Brill-Noether, Theorem \ref{Brill-Noether theorem}]
Let $X_m$ be a del Pezzo surface with $m \leq 5$. Let $\bv$ be a Chern character such that $\nu(\bv).L \geq -2$ and $\cP_{L^W}(\bv) \neq \emptyset$.
\begin{enumerate}
\item[(1)] If $\nu(\bv).C \geq -1$ for all $(-1)$-curves $C$, then $\bv$ is non-special.
\item[(2)] Let $W_m$ be the Weyl group acting on $\Pic(X_m)$. If there exists some $\sigma \in W_m$ such that $\nu(\bv).\sigma(L) \leq -1$ or $\nu(\bv).\sigma(L-E_i) \leq -1$ for some $i$, then $\bv$ is non-special.
\item[(3)] Let $\nu(\bv).\sigma(L) > -1$ and $\nu(\bv).\sigma(L-E_1) > -1$ for all $\sigma \in W_m$. Suppose $C$ is a $(-1)$-curve such that $\nu(\bv).C < -1$ and let $\pi$ be the map contracting $C$. Then $\bv$ is non-special if and only if $\pi_*(\cE)$ is non-special for a general sheaf $\cE$ and $\chi(\pi_*(\cE)) \leq 0$. 
\end{enumerate}
\end{theorem}

We give a method for constructing direct sums of line bundles of any given total slope $\nu \in \Pic(X_m)_{\QQ}$. We show that when $\nu(\bv)$ is nef, the summands in our construction are close to being nef, so that their higher cohomology vanishes by a theorem of Castelnuovo. To obtain non-special sheaves of the same total slope and larger discriminant, we use elementary modifications. When $\nu(\bv)$ is not nef, we blowdown $(-1)$-curves $C$ such that $\nu(\bv).C < 0$ and use the Leray spectral sequence to compute cohomology.

On the blowup of $\PP^2$ at a general collection of points, the cohomology of line bundles has been well-studied. This is the subject of the Segre-Harbourne-Gimigliano-Hirchowitz (SHGH) conjecture. The conjecture is known for at most $9$ points by work of Castelnuovo, but it is open in general. Our results provide evidence for the following higher rank formulation of the SHGH conjecture.

\begin{conj}[Conjecture 1.7 \cite{CoskunHuizengaWeakBrillNoether}]
\label{Higher rank SHGH}
Let $X_m$ be the blowup of $\PP^2$ at $m$ general points, and let $H$ be an ample class such that $H.(K_{X_m}+F)< 0$, with $F$ a fiber class. If $M_{X,H}(\bv) \neq \emptyset$ and $\ch_1(\bv)$ is nef, then $\bv$ is non-special.
\end{conj}

\subsection{Existence}
In \cite{DrezetLePotierExistenceP2} for $\PP^2$, Dr\'ezet and Le Potier show that the existence problem is controlled by the exceptional bundles on $\PP^2$. They construct a function $\delta: \RR \to \RR$ whose graph in the $(\mu , \Delta)$-plane completely determines when $M(\bv)$ is nonempty. If $(\mu(\bv),\Delta(\bv))$ lies above the graph of $\delta$, then $M(\bv)$ is nonempty. Otherwise, $M(\bv)$ is empty or $\bv$ is the Chern character of a semi-exceptional bundle. For Hirzebruch surfaces and generic polarizations, the classification of nonempty moduli spaces is known by \cite{Hirzebruchexist}. 
Rudakov gave existence theorems for del Pezzo surfaces with the anti-canonical polarization in \cite{RudakovQuadricSemistable, RudakovExistenceDelPezzo} conditional on the existence of restricted smooth complete families. He shows their existence for $X_1$ in \cite{RudakovExistenceDelPezzo}. As a consequence of our construction in Section \ref{good bundle construction and properties}, we are able to construct restricted smooth complete families for del Pezzo surfaces of degree at least $3$. Before stating the result, we 
 recall the definition of \textit{Dr\'ezet  - Le Potier condition} 

A torsion-free coherent sheaf $\cF$ (or Chern character) satisfies the \textit{Dr\'ezet  - Le Potier condition} (abbr. as \textit{DL condition}) if 

\begin{enumerate}[label=DL\arabic*]
\item  For every exceptional bundle $\cE$ satisfying $r(\cE) < r(\cF)$ and
\[
\mu(\cF) \leq \mu(\cE) \leq \mu(\cF) + K_{X_m}^2,
\]
we have $\chi(\cE,\cF) \leq 0$, \\

\item  and  for every exceptional bundle $\cE$ satisfying $r(\cE) < r(\cF)$ and
\[
\mu(\cF)-K_{X_m}^2 \leq \mu(\cE) \leq \mu(\cF),
\]
we have $\chi(\cF,\cE) \leq 0$.
\end{enumerate}

\begin{theorem}[Theorem \ref{existence del Pezzo}]
Let $X_m$ be a del Pezzo surface with $m \leq 6$. Let $\bv \in K(X_m)$. If $\Delta(\bv)>\frac{1}{2}$, then there exists a $-K_{X_m}$-stable bundle of Chern character $\bv$ if and only if $\bv$ satisfies the DL condition. 
\end{theorem}

\subsection{Outline} The outline of the paper is as follows. In Section \ref{prelims}, we review some basic facts about stability conditions and moduli stacks of sheaves. Section \ref{del Pezzo stuff} covers some basic cohomology computations for line bundles on del Pezzo surfaces and includes a brief discussion of the SHGH Conjecture.

Section \ref{prioritary section} compiles facts about prioritary sheaves, some of which can be found in \cite{WalterComponentsofStack} and \cite{WalterIrreducibility}. These are the main classes of sheaves we work with throughout the paper. We discuss how the numerical invariants and cohomology of prioritary sheaves behave under blowups and elementary modifications and end with the proof of Proposition \ref{refined Bogomolov}.

Section \ref{good bundle construction and properties} deals with the construction of good bundles and their numerical properties. The classification result is proved in Section \ref{Brill-Noether}. We compute the cohomology of the bundles constructed in Section \ref{good bundle construction and properties} after twisting by line bundles. The results obtained in Section \ref{prioritary section} take care of the remaining cases with higher discriminant.

We briefly recall the main definitions and results from \cite{RudakovExistenceDelPezzo} in Section \ref{DL condition and Rudakov theorems}. The construction of restricted smooth complete families and the proof of the existence of stable and semistable sheaves occupy all of Section \ref{restricted sheaves and existence}.

\subsection{Acknowledgements} The first author would like to thank Jack Huizenga for his support, encouragement, and helpful discussions throughout the project. The second author is grateful to Penn State University for their hospitality and inspiring atmosphere. The authors would also like to thank Arend Bayer, Alexey Gorodentsev, and Dmitrii Pedchenko for helpful discussions.\\

\section{Preliminaries}
\label{prelims}

The following is a brief overview of our notation and stability conditions. The standard reference for moduli of sheaves is \cite{HuybrechtsLehnModuliOfSheaves}. Our remarks in the case of surfaces follows mostly \cite{Huizenga-survey}.

\subsection*{Conventions} All varieties considered in this paper are smooth projective over $\CC$, and all sheaves are coherent unless specified otherwise. For a variety $X$ and coherent sheaves $\cE$ and $\cF$, we set $h^i(X,\cE) = \dim H^i(X,\cE)$, $\hom(\cE,\cF) = \dim \Hom(\cE,\cF)$, and $\ext^i(\cE,\cF) = \dim \Ext^i(\cE,\cF)$.

\subsection{Numerical invariants}

Let $X$ be a smooth projective surface. We denote by $K(X)_{\QQ}$ the Grothendieck group of $X$ modulo numerical equivalence with coefficients in $\QQ$. Let $\cE$ be a torsion-free sheaf on $X$. Then $\ch(\cE)=(\ch_0(\cE),\ch_1(\cE),\ch_2(\cE))$, where
\begin{align*}
\ch_0(\cE) & = r(\cE) \hspace{2mm}\text{(rank of $\cE$)}\\
\ch_1(\cE) & = c_1(\cE)\\
\ch_2(\cE) & = \frac{c_1^2-2c_2}{2},
\end{align*}
and where $c_i$ are the usual Chern classes. It is useful to repackage $\ch(\cE)$ as logarithmic invariants of the form
\[
\ch(\cE) = (r(\cE),\nu(\cE),\Delta(\cE)),
\]
where $\nu(\cE) = \ch_1(\cE)/\ch_0(\cE)$ is the \textit{total slope of $\bv$} and
\[
\Delta(\cE) = \frac{\nu(\cE)^2}{2} - \frac{\ch_2(\cE)}{r(\cE)}
\]
is the \textit{discriminant of $\bv$}. For any torsion-free sheaves $\cE$ and $\cF$, they satisfy the properties
\begin{align*}
r(\cE \otimes \cF) & = r(\cE)r(\cF),\\
\nu(\cE \otimes \cF) & = \nu(\cE) + \nu(\cF),\\
\text{and} \hspace{1cm} \Delta(\cE \otimes \cF) & = \Delta(\cE) + \Delta(\cF).
\end{align*}
It is important to note that when $\cL$ is a line bundle, we have $\Delta(\cL) = 0$ and $\Delta(\cE \otimes \cL) = \Delta(\cE)$.

We can formally generalize Riemann-Roch for line bundles on surfaces to divisors with rational coefficients by defining
\[
P(\nu):= \chi(\cO_X) +\frac{\nu.(\nu-K_X)}{2}.
\]
Let $\chi(\cE, \cF) = \sum (-1)^i \ext^i(\cE, \cF)$. We have the following Riemann-Roch formula for torsion-free sheaves $\cE$ and $\cF$:
\[
\chi(\cE,\cF) = r(\cE)r(\cF)(P(\nu(\cF)-\nu(\cE))-\Delta(\cE)-\Delta(\cF)).
\]
When $\cE = \cO$, this reduces to the formula
\[
\chi(\cF) = r(\cF)(P(\nu(\cF)) - \Delta(\cF)).
\]

\subsection{Stability}
\label{stability}
Let $\cE$ be a torsion-free sheaf on a polarized smooth projective variety $(X,H)$. We can consider the \textit{reduced Hilbert polynomial}
\[
p_{\cE}(m) := \frac{\chi(\cE(m))}{r(\cE)}.
\]
We say $\cE$ is \textit{$H$-semistable} (or just \textit{semistable} when $H$ is implicitly understood) if $p_{\cF}(m) \leq p_{\cE}(m)$ for every nonzero subsheaf $\cF$ of strictly smaller rank and for all $m \gg 0$. If the inequality is strict for all such $\cF$, we say $\cE$ is \textit{$H$-stable} (or just \textit{stable}). Every strictly $H$-semistable sheaf $\cE$ has a Jordan-H\"older filtration
\[
0 \subset \cE_1 \subset \cdots \subset \cE_k = \cE,
\]
with stable factors $\gr_{JH}^i(\cE) : = \cE_i/\cE_{i-1}$, which are unique up to isomorphism and reordering of factors. Two $H$-semistable sheaves are \textit{$S$-equivalent} if they have the same Jordan-H\"older factors up to isomorphism and reordering.

Fix $\bv \in K(X)$. There is a stack $\cM_{X,H}(\bv)$ called the \textit{moduli stack of $H$-semistable sheaves} of Chern character $\bv$, which is an open algebraic substack of the stack $\mathscr{C}oh_X(\bv)$ (see Theorem 4.6.2.1 in \cite{LaumonMoretBaillyStacks}), by openness of semistability. It is corepresented by a projective scheme $M_{X,H}(\bv)$, called the \textit{moduli space of $H$-semistable sheaves}, which parametrizes sheaves up to $S$-equivalence. Although our results apply to $M_{X,H}(\bv)$, we will work primarily on the level of stacks, where we can avoid technical issues with purely semistable moduli spaces.

We will also work with a simpler, but related notion, of stability. Let $\cE$ be a torsion-free sheaf on $(X,H)$. Define the \textit{$\mu_H$-slope} of $\cE$ to be 
\[
\mu_H(\cE):= \frac{\nu(\cE).H}{H^2}.
\]
We say $\cE$ is \textit{$\mu_H$-semistable} (or just \textit{slope semistable} when $H$ is implicitly understood) if for every nonzero subsheaf $\cF$ of strictly smaller rank, we have $\mu_H(\cF)\leq \mu_H(\cE)$. If the inequality is always strict, then $\cE$ is \textit{$\mu_H$-stable} (or \textit{slope stable}). Slope semistability and slope stability are also open conditions.

We can rewrite $p_{\cE}(m)$ in the following way:
\[
p_{\cE}(m) = \frac{1}{2}m^2 + \frac{H.(\nu(\cE) - \frac{1}{2}K_X)}{H^2}m + \frac{\chi(\cE)}{r(\cE)H^2}.
\]
From this, it is clear we have the implications
\[
\mu_H\text{-stable} \Rightarrow \text{stable} \Rightarrow \text{semistable} \Rightarrow \mu_H\text{-semistable}.
\]
Let $\cM^{\mu ss}_{X, H}(\bv)$ be the moduli stack of $\mu_H$-semistable sheaves of Chern character $\bv$. Then $\cM_{X,H}(\bv)$ is an open substack of $\cM^{\mu ss}_{X, H}(\bv)$. It is important to note that for any line bundle $\cL$, the sheaf $\cE \otimes \cL$ is $\mu_H$-(semi)stable if and only if $\cE$ is. This is not true in general for $H$-(semi)stable sheaves.

We write $\nu_{H,\max}(\cE)$ ($\mu_{H,\max}(\cE)$) and $\nu_{H,\min}(\cE)$ ($\mu_{H,\min}(\cE)$) for the total slopes ($H$-slopes) of the Gieseker Harder-Narasimhan factors of maximal and minimal reduced Hilbert polynomial, respectively.\\

\subsection{Stacks} All flat families of sheaves we consider live inside the algebraic stack $\mathscr{C}oh_X$, the stack of all coherent sheaves on $X$. Let $\cF$ be a flat family of sheaves on $X$ parametrized by a smooth variety $S$. For each point $p \in S$, there is a \textit{Kodaira-Spencer map} (see 10.1.8 \cite{HuybrechtsLehnModuliOfSheaves})
\[
\omega_p:T_pS \to \Ext^1(\cF_p,\cF_p).
\]
A family of sheaves is \textit{complete} if the Kodaira-Spencer map surjects for each point $p\in S$. In the setting of moduli theory, it is often the case that the family of all semistable sheaves is unwieldy, while complete families can be explicitly described in terms of resolutions or extensions. However, they have the disadvantage of not always containing all isomorphism classes of semistable sheaves because of their local nature. On the other hand, the stacks we consider are all smooth and open in $\mathscr{C}oh_X$. Thus, for the purposes of computing the cohomology of the generic sheaf, it is enough to work with complete families.

\section{del Pezzo surfaces}
\label{del Pezzo stuff}

The proof of Brill-Noether relies heavily on the knowledge of the cohomology of line bundles, which is known in the case of del Pezzo surfaces. We review some definitions for del Pezzo surfaces and discuss some of their properties.

A del Pezzo surface is a Fano variety of dimension $2$, i.e. a smooth projective surface with ample anti-canonical bundle. Let $X_m$ be the blowup of $\PP^2$ at $m$ general points. Up to isomorphism, del Pezzo surfaces are either $X_m$ for $m<9$ or $\PP^1 \times \PP^1$. For blowups, we have
\begin{gather*}
\Pic(X_m) = \ZZ L \oplus \ZZ E_1 \oplus \cdots \oplus \ZZ E_m, \\
\text{and }L^2 = 1, \hspace{1mm} E_i^2 = -1, \hspace{1mm} L.E_i = E_i.E_j = 0 \text{ for } j \neq i.
\end{gather*}
The divisor $L$ denotes the pullback of the line class on $\PP^2$ along the blowup map, and the $E_i$ represent the exceptional divisors obtained from blowing up the points $p_i$ on $\PP^2$. For $\PP^1 \times \PP^1$, the Picard group is generated by the two fiber classes, denoted by $F_1$ and $F_2$, with $F_1^2 = F_2^2 = 0$ and $F_1.F_2 = 1$.

The canonical divisor for $X_m$ is $K_{X_m} = -3L + \sum_{i=1}^m E_i$. If $D = aL - \sum_{i=1}^m b_iE_i$, then 
\[
\chi(\cO_X(D)) = \frac{(a+1)(a+2)}{2} - \sum_{i=1}^m \frac{b_i(b_i+1)}{2}
\]
by Riemann-Roch. In the case of the quadric $\PP^1 \times \PP^1$, the Riemann-Roch theorem for a line bundle $\cO(a,b)$, where $a$ and $b$ denote the coefficients of $F_1$ and $F_2$, respectively, takes on the form
\[
\chi(\cO(a,b)) = (a+1)(b+1).
\]

A $(-1)$-curve on a surface is a smooth rational curve whose self-intersection is  $-1$. There are only finitely many $(-1)$-curves on del Pezzo surfaces, and they carry certain symmetries. We recall their properties here.

\begin{prop}[\cite{ManinWeylGroup}]
\label{Weyl group properties}
For each integer $m$ such that $2 \leq m\leq 9$, there is a root system on $\Pic(X_m)_{\RR}$ and a corresponding Weyl group $W_m$ that acts on $\Pic(X_m)_{\RR}$ with the following properties.
\begin{enumerate}
\item[(i)] The intersection pairing is preserved. That is, for all $\sigma \in W_m$ and divisors $D_1$ and $D_2$ we have $\sigma(D_1).\sigma(D_2) = D_1.D_2$.
\item[(ii)] $W_m$ fixes the canonical class $K_{X_m}$.
\item[(iii)] For each integer $n$ such that $1 \leq n \leq m$ and $n \neq m-1$, the group $W_m$ acts transitively on the set of collections of pairwise disjoint $(-1)$-curves of size $n$. When $n=m-1$, there are two orbits, one consisting of collections that can be extended to lists of size $m$ and the other consisting of collections that are maximal in size.
\end{enumerate}
\end{prop}

We note that the two orbits appearing in Proposition \ref{Weyl group properties} (iii) for $n=m-1$ correspond to the exceptional divisors belonging to the blowups of either $\PP^1 \times \PP^1$ or $X_1$. When $m \geq 2$, the pseudo-effective cone $\overline{\Eff}(X_m)$ is generated by the $(-1)$-curves of $X_m$. If $m=1$, we have $\overline{\Eff}(X_1)=\langle E_1,L-E_1\rangle$. We list the $(-1)$-curves up to $X_6$ up to interchanging the $E_i$ by the Weyl group action.
\begin{enumerate}
\item[(1)] $E_1$
\item[(2)] $L-E_1-E_2$
\item[(3)] $2L-\sum_{i=1}^5 E_i$
\end{enumerate}

By cone duality, the nef divisors on $X_m$ are precisely the ones nonnegative on $\overline{\Eff}(X_m)$. Since $\overline{\Eff}(X_m)$ is finitely generated, we will often show that a divisor $D$ is nef by showing $D.C \geq 0$ for every generator $C$ of $\overline{\Eff}(X_m)$, which, for $m \geq 2$, amounts to showing that $D.C \geq 0$ for all $(-1)$-curves $C$. By Proposition \ref{Weyl group properties} (i), the Weyl group action preserves $\overline{\Eff}(X_m)$ and the nef cone $\Nef(X_m)$. If $m \geq 2$, the primitive integral divisors on the extremal rays of $\Nef(X_m)$ split into two orbits under the Weyl group action. One orbit consists of the translates of the class $L$, and the other orbit consists of the translates of the fiber class $L-E_1$.

\begin{prop}
\label{h0h2vanishing}
Let $D = aL - \sum_{i=1}^m b_iE_i$.
\begin{enumerate}
\item[(a)] If $a < 0$, then $H^0(X_m,\cO_{X_m}(D)) = 0$.\\
\item[(b)] If $a > -3$, then $H^2(X_m,\cO_{X_m}(D)) = 0$.
\end{enumerate}
\end{prop}

\begin{proof}
Note that $L$ is nef and $L.D = a$. If $a < 0$, then $D$ is not effective. By Serre duality, this implies (b) as well.
\end{proof}

\begin{prop}
\label{higherdirectimages}
Let $f:S \to S'$ be the blowup of a smooth projective surface $S'$ at a point with exceptional divisor $E$.
\begin{enumerate}
\item[(a)] $R^1f_*\cO_S(kE) = 0$ when $k \leq 1$, and it is nonzero otherwise.\\
\item[(b)] $f_*\cO_S(kE) = \cO_{S'}$ when $k \geq 0$.
\end{enumerate}
\end{prop}
\begin{proof}
The proof for (b) follows from the pushforward of the restriction sequence
\[
0 \to \cO_S((k-1)E) \to \cO_S(kE) \to \cO_E(-k) \to 0 
\]
and induction on $k$.

For (a), we adapt the proof for the case $k=0$ given in Proposition V 3.4 in \cite{HartshorneAlgebraicGeometry}. Note that $R^1f_*\cO(kE)$ is a coherent sheaf supported at the center of the blowup. By the theorem on formal functions, it vanishes if and only if $\varprojlim H^1(E_n,\cO_{E_n}(kE)) = 0$, where $E_n$ is the closed subscheme of $S$ defined by $\cO(-nE)$. For each $n \geq 1$, there are exact sequences (on $S$)
\[
0 \to \cO_E(n) \to \cO_{E_{n+1}} \to \cO_{E_n} \to 0.
\]
Tensoring with $\cO(kE)$ and induction on $n$ gives the result.
\end{proof}

The most difficult cohomology computations involve nef line bundles (or, more generally, effective linear series that are not ($-1$)-special). This is the subject of the famous Segre-Harbourne-Gimigliano-Hirschowitz (SHGH) conjecture and its many variations. See \cite{CHMRVariationsOfNagata} and \cite{CilibertoMirandaDegenerations} for more details and progress on the problem.

\begin{mydef}
A line bundle $\cO(D)$ on $X_m$ is called $(-1)$\textit{-special} if $D$ is effective and $D.C \leq -2$ for some $(-1)$-curve $C$.
\end{mydef}

Geometrically, $(-1)$-special line bundles correspond to effective divisors with nonreduced base loci.

\begin{conj}[SHGH Conjecture]
If $D$ is effective, then $\chi(\cO(D)) = \dim(H^0(X_m,\cO(D)))$ if and only if $\cO(D)$ is not $(-1)$-special.
\end{conj}
Note that the vanishing of $H^2(X_m,\cO(D))$ is a consequence of Proposition \ref{h0h2vanishing} and positivity conditions on nef divisors, so the difficulty in the statement lies in showing that $H^1(X_m,\cO(D))=0$. For del Pezzo surfaces and the blowup at 9 general points, the full SHGH conjecture is known and has been proved first by Castelnuovo. More recent proofs can be found in \cite{NagataSHGH} and \cite{HarbourneSHGH}.

\begin{theorem}
\label{SHGHDelPezzo}
The SHGH conjecture is true for del Pezzo surfaces and the blowup at $9$ general points. \qed
\end{theorem}

\begin{corollary}
\label{line bundle cohomology}
Let $2 \leq m \leq 9$, and suppose $D = aL -\sum_{i=1}^{m} b_iE_i$ with $a \geq -2$. Suppose $D.C \geq -1$ for all $(-1)$-curves $C$. If the set of $(-1)$-curves $C$ with $C.D=-1$ are pairwise disjoint, then $H^i(X_m,\cO(D)) = 0$ for all $i > 0$.
\end{corollary}

\begin{proof}
The assumption on $a$ and Proposition \ref{h0h2vanishing} shows $H^2(X_m,\cO(D)) = 0$. For any $C$ such that $D.C = -1$, there is an exact sequence
\[
0 \to \cO(D-C) \to \cO(D) \to \cO_C(-1) \to 0.
\]
Let $f:X_m \to X_{m-1}$ be the map contracting $C$, and apply the functor $f_*$ to the exact sequence. We deduce that $R^if_*\cO(D) \cong R^if_*\cO(D-C)$ for all $i$. Since $(D-C).C=0$, the line bundle $\cO(D-C)$ is the pullback of a line bundle on $X_{m-1}$, and we obtain $R^1f_*\cO(D-C) = 0$. Thus, the cohomology of $\cO(D)$ is the cohomology of $f_*\cO(D)$.

Let $\pi: X_m \to Y$ be the map contracting the set of disjoint $(-1)$-curves $C$ where $D.C = -1$. If $Y = X_n$ and $n \geq 2$, then $\pi_*D$ is nef, so $\pi_*\cO(D)$ has no higher cohomology by Theorem \ref{SHGHDelPezzo}. If $n < 1$ or $Y = \PP^1 \times \PP^1$, then we can reduce to the case $m = 2$ and split into two subcases.
\begin{enumerate}
\item[1)] If $D.E_i = -1$ for at least one $i$, then $D.(L-E_1-E_2) \geq 0$. Without loss of generality, we can assume $D.E_2 = -1$. Then $D.(L-E_1) = D.(L-E_1-E_2 + E_2) \geq -1$. Thus, taking $Y = X_1$, we see that $\chi(\pi_*\cO(D)) \geq 0$, so the higher cohomology vanishes.
\item[2)] If $D.(L-E_1-E_2) = -1$, then $D.E_i \geq 0$ for both $i$. Thus, taking $Y = \PP^1 \times \PP^1$, we get $\pi_*\cO(D) \cong \cO_Y(a,b)$ with $a \geq -1$ and $b \geq -1$, so the higher cohomology vanishes.
\end{enumerate}
\end{proof}

\section{Prioritary sheaves}
\label{prioritary section}

\subsection{Prioritary sheaves and Walter's theorem}
\begin{mydef}
\label{prioritarydef}
Let $X$ be a surface and let $D$ be a divisor on $X$. A torsion-free sheaf $\cE$ is \textit{$D$-prioritary} if $\Ext^2(\cE,\cE(-D)) = 0$.
\end{mydef}

\begin{rmk}
\label{prioritary line bundle twist}
If $\cL$ is any line bundle, then $\cE \otimes \cL$ is $D$-prioritary if and only if $\cE$ is $D$-prioritary.
\end{rmk}

Since being torsion-free and $D$-prioritary are open conditions, there is an open algebraic substack $\cP_{X,D}(\bv)$ ($\cP_D(\bv)$ if the surface is unambiguous) of $\mathscr{C}oh_{X}(\bv)$ consisting of $D$-prioritary sheaves. The advantages of using prioritary sheaves are twofold. First, they are easier to construct and have good functorial properties. Second, we have the following lemma:

\begin{lemma}[Lemma 4 \cite{WalterIrreducibility}]
\label{prioritary smooth restriction morphism}
Let $D$ be an effective divisor on a surface $X$. Let $i:D \hookrightarrow X$ be a closed immersion. If $\cE$ is $D$-prioritary and $i^*\cE$ is locally free, then the morphism of stacks $i^*: \cP_{D}(\bv) \to \mathscr{C}oh_{D}(i^*\bv)$ defined by $\cF \mapsto i^*\cF$ is smooth near $\cE$. \qed
\end{lemma}

A surface $X$ is \textit{geometrically ruled} if there exists a smooth morphism $\pi:X \to C$ to a curve $C$ with all fibers $\PP^1$. If $X$ is birational to a geometrically ruled surface $X'$, then $X$ is \textit{birationally ruled}. It is a standard fact that a birationally ruled surface is either $\PP^2$, geometrically ruled, or a blowup of a geometrically ruled surface. Thus, if $X\neq \PP^2$, there is a flat surjective map $\pi:X \to C$ to a curve $C$ with general fiber $\PP^1$.

Let $p \in C$, and let $F$ be the class of $\pi^{-1}(p)$ in $\NS(X)$. If $H.(K_X+F) < 0$ and $\cE$ is $\mu_H$-semistable, then $\cE$ is $F$-prioritary. Indeed, we have $\Ext^2(\cE,\cE(-F)) \cong \Hom(\cE,\cE(K_X+F))^*$ by Serre duality, and $\Hom(\cE,\cE(K_X+F)) = 0$ by stability. We need the following theorem in order to speak of general sheaves in moduli spaces.

\begin{theorem}[Theorem 1, Proposition 2 \cite{WalterIrreducibility}]
\label{irreducibility}
Let $\pi:X \to C$  be a birationally ruled surface, and let $\bv \in K(X)$ with $\ch_0(\bv) \geq 2$. Then $\cP_F(\bv)$ is smooth and irreducible. Moreover, if $H$ is a polarization such that $H.(K_X+F) < 0$ and $M_{X,H}(\bv)$ is nonempty, then $M_{X,H}(\bv)$ is irreducible and normal.
\qed
\end{theorem}

\begin{rmk}
\label{other stacks irreducible}
Let $(X,H)$ be as in the statement of Theorem \ref{irreducibility}. If $\cM_{X,H}^{\mu ss}(\bv)$ and $\cM_{X,H}(\bv)$ are nonempty, then they are irreducible by Theorem \ref{irreducibility}.
\end{rmk}

\begin{prop}
\label{effective+prioritary}
Let $D_1$ and $D_2$ be divisors on a surface $X$. Let $\cE$ be a $D_1$-prioritary sheaf. If $D_1-D_2$ is effective, then $\cE$ is $D_2$-prioritary.
\end{prop}

\begin{proof}
Let $\cE$ be a torsion-free sheaf on $X$, and let $Y$ be a codimension $1$ subscheme. Consider the restriction sequence
\[
0 \to I_Y \to \cO_X \to \cO_{Y} \to 0.
\]
After tensoring the sequence with $\cE$, we get an injection $\mathcal{T}or^{\cO_X}_1(\cE,\cO_{Y}) \to \cE \otimes I_Y $ from a torsion sheaf into a torsion-free sheaf, so $\mathcal{T}or^{\cO_X}_1(\cE,\cO_{Y}) = 0$. Now, apply the functor $\Hom(\cE(D_2),-)$ to the sequence 
\[
0 \to \cE(-(D_1-D_2)) \to \cE \to \cE|_{D_1-D_2} \to 0.
\]
The conclusion is evident from the long exact sequence
\[
\cdots \to \Ext^2(\cE,\cE(-D_1)) \to \Ext^2(\cE,\cE(-D_2)) \to 0.
\]
\end{proof}

\begin{rmk}
\label{semistableimpliesprior}
The choice of a ruling on $X_m$ is equivalent to a choice of a map coming from a divisor $F$ such that $F^2=0$ and the complete linear series is isomorphic to $\PP^1$ and contains a smooth rational curve. We call such a divisor a \textit{fiber class} of $X$.

Let $m \geq 1$. The class $F_i=L-E_i$ is easily seen to be a fiber class. Since $L - F_i = E_i$ is effective, we have an inclusion of open substacks $\cP_L(\bv) \subset \cP_{F_i}(\bv)$ by Proposition \ref{effective+prioritary}. When $m \leq 6$, the divisor $K_{X_m}+F_i$ is anti-effective, so the hypotheses of Theorem \ref{irreducibility} are satisfied for every ample $H$. Thus, the moduli space $M_{H}(\bv)$ is irreducible for any polarization $H$. If $m = 7$, there exist polarizations $H$ such that $H.(K_{X_7}+F) \geq 0$. For example, the anti-canonical polarization gives $-K_{X_7}.(K_{X_7}+F) = 0$, and when $n \geq 3$ the divisor $-nK_{X_7}-F$ is ample and $(-nK_{X_7}-F).(K_{X_7}+F) > 0$.

When $m\leq 5$, we have $K_{X_m}+L$ is anti-effective, so $\mu_H$-semistable sheaves are $L$-prioritary. We will be primarily concerned with the stack $\cP_L(\bv)$. Note in Theorem \ref{irreducibility} for $X = X_m$, the statement remains true if we replace $F$ with $L$.
\end{rmk}

\subsection{Bogomolov-type statement}

\begin{mydef}
Let $C$ be a curve on a surface $X$. Let $\cE$ be a sheaf on $X$ locally free in a neighborhood of $C$, and let $\cF$ be a locally free sheaf on $C$. We call a sheaf $\cE'$ a \textit{type $1$ elementary transformation of $\cE$ along $\cF$} if $\cE'$ is the kernel of a surjective map $\cE \to \cF$. If $\cE'$ is locally free in a neighborhood of $C$ and $\cE'$ is given by an extension in $\Ext^1(\cF, \cE)$, then we say $\cE'$ is a \textit{type $2$ elementary transformation of $\cE$ along $\cF$}.
\end{mydef}

In the literature, elementary transformations always refer to our type $1$ elementary transformations. We define the notion of a type $2$ elementary transformation for convenience, since we will usually work with bundles given by an extension. Note in the type $1$ case, the sheaf $\cE'$ is a locally free near $C$ by a local computation.

The remainder of this section is dedicated to proving the following Bogomolov-type statement.

\begin{prop}
\label{refined Bogomolov}
Let $\cV$ be a sheaf on $X_m$ of rank $r \geq 2$. Suppose $\cV$ is a type $2$ elementary transformation of $\cO_{X_m}(-2L)^{\oplus a} \oplus \cO_{X_m}(-L)^{\oplus b}$ along $\bigoplus_{i=1}^m \cO_{E_i}(-1)^{\oplus d_i}$, where $m \geq 0$, $a,b \geq 0$, and $0 \leq d_i < r$. Then $\cV$ is $L$-prioritary, and for any $D \in \Pic(X_m)$, there are no $L$-prioritary sheaves of the same rank and total slope as $\cV \otimes \cO(D)$ with strictly smaller discriminant. In particular, for $m \leq 5$ and any polarization $H$, there are no $\mu_H$-semistable sheaves of the same rank and total slope as $\cV \otimes \cO(D)$ with strictly smaller discriminant.
\end{prop}

\begin{rmk}
\label{all total slopes possible}
For any rank $r \geq 1$ and total slope $\nu \in \Pic(X_m)_{\QQ}$, it is not hard to see that there is a divisor $D$ and a type $2$ elementary transformation $\cV$ of $\cO_{X_m}(-2L)^{\oplus a} \oplus \cO_{X_m}(-L)^{\oplus b}$ along $\bigoplus_{i=1}^m \cO_{E_i}(-1)^{\oplus d_i}$ for some integers $a$, $b$, $d_1$, \ldots, and $d_m$ such that $r(\cV(D)) = r $ and $\nu(\cV(D)) = \nu$. Furthermore, the divisor $D$ and integers $a$, $b$, $d_1$, \ldots, and $d_m$ are unique even though there is a choice for the type $2$ elementary transformation $\cV$. Since the Chern character of $\cV$ is determined by the integers $a$, $b$, $d_1$, \ldots, and $d_m$, there is a well-defined function
\begin{align*}
\Delta_m: \NN \times \Pic(X_m)_{\QQ} \to \QQ \\
\Delta_m(r,\nu)  = \Delta(\cV(D)).
\end{align*}
If $m \leq 5$ and $\cE$ is a $\mu_H$-semistable bundle on $X_m$, then $\Delta(\cE) \geq \Delta_m(r(\cE),\nu(\cE))$ by Proposition \ref{refined Bogomolov}.
\end{rmk}

There are two key ingredients to the proof. First, we will show that $\cV$ has minimal discriminant among $L$-prioritary sheaves of the same rank and total slope. Then we will give a formula for relating the discriminant of a sheaf $\cE$ and an elementary transformation of $\cE$.

Note that the bundle $\cO_{X_m}(-2L)^{\oplus a} \oplus \cO_{X_m}(-L)^{\oplus b}$ is $L$-prioritary by direct computation.

\begin{prop}
\label{No L-prioritary sheaves with smaller discriminant}
Let $\cV = \cO^{\oplus a} \oplus \cO(L)^{\oplus b}$ be a vector bundle on $\PP^2$, and let $n$ be an integer. Suppose $\cW$ is a sheaf such that $r(\cW) = r(\cV(nL))$ and $\nu(\cW) = \nu(\cV(nL))$. If $\Delta(\cW) < \Delta(\cV(nL))$, then $\cW$ is not $L$-prioritary.
\end{prop}

\begin{proof}
It is enough to prove the statement when $n=0$. Every $L$-prioritary sheaf on $\PP^2$ has a nonpositive Hilbert polynomial by Proposition 1.1(3) \cite{GottscheHirschowitzBrillNoetherP2}. The minimal value of the Hilbert polynomial $P_{\cV}(m)$ is $P_{\cV}(-2) = 0$, and $P_{\cV}(-1) = 0$ (resp. $P_{\cV}(-3) = 0$) when $b =0$ (resp. $a=0$). Any sheaf of smaller discriminant with the same rank and total slope as $\cV$ has a strictly positive Hilbert polynomial, thus is not $L$-prioritary.
\end{proof}

An elementary computation shows the bundle $\cV$ in Proposition \ref{No L-prioritary sheaves with smaller discriminant} has no infinitesimal deformations, i.e. $\Ext^1(\cV,\cV) = 0$. Thus, if $\ch(\bv) = \ch(\cV)$, then the general object in $\cP_L(\bv)$ is isomorphic to $V$.

We now proceed with the second key component.

\begin{prop}[Lemma 6 \cite{WalterIrreducibility}]
\label{elementarytransformationexact}
Let $f: X' \to X$ be the blowup of a surface $X$ at a point $x \in X$. Let $E$ be the exceptional divisor in $X'$. Suppose that $\cE$ is a sheaf of rank $r$ on $X'$ such that $\cE|_E \cong \cO_E^{r-d} \oplus \cO_E(-1)^d$ for some $d$. Then $f_*(\cE)$ is locally free in a neighborhood of $x$, and there are exact sequences
\begin{gather*}
0 \to f^*f_*(\cE) \to \cE \to \cO_E(-1)^{\oplus d} \to 0,\\
0 \to \cE(-E) \to f^*f_*(\cE) \to \cO_E^{\oplus r-d} \to 0.
\end{gather*}
Moreover, for any divisor $D$ on $X$, we have $\Ext_{X'}^2(\cE,\cE(f^*(D))) \cong \Ext_{X}^2(f_*(\cE),f_*(\cE)(D))$.
\qed
\end{prop}

\begin{rmk}
\label{Grassmannian bundle over stack}
Let $X$ be a birationally ruled surface, and let $f$, $X'$, $x$, $E$, and $\cE$ be as in Proposition \ref{elementarytransformationexact}. Denote by $\bv$ and $\bv'$ the Chern characters of $\cE$ and $f_*(\cE)$, respectively, and assume $\ch_0(\bv) \geq 2$. If $\cP_F(\bv')$ is nonempty, then the substack $\cP_{F,0}(\bv') \subset \cP_F(\bv')$ of locally free sheaves is open dense by \cite{WalterIrreducibility}. More generally, if $\cP_F(\bv')$ is nonempty, there is an open dense substack $\cP_{F,(x,0)}(\bv')$ of sheaves locally free at $x$.

On the other hand, we have $\cP_{F}(\bv) \subset \cP_E(\bv)$ by Proposition \ref{effective+prioritary}. Now let $C$ be a smooth rational curve on a surface $X$. If $\cP_F(\bv)$ is nonempty, then there is an open dense substack $\cR_F^E(\bv)$ of sheaves of rigid splitting type along $E$ (the restriction of a sheaf $\cE$ to $E$ is $\cE|_E \cong \cO_E(n)^{\oplus a} \oplus \cO_E(n+1)^{\oplus b}$ for some integer $n$ and nonnegative integers $a$ and $b$) by Lemma \ref{prioritary smooth restriction morphism}. Thus, Proposition \ref{elementarytransformationexact} implies that the stack $\cR_F^E(\bv)$ is a $G(d,r)$-bundle over $\cP_{F,(x,0)}(\bv')$. If $\cF \in \cP_{F,(x,0)}(\bv')$, the choice of a $d$-plane in $H^0(E,\cF|_E)$ uniquely determines the isomorphism class of the extension. Moreover, the stack $\cP_F(\bv)$ is nonempty if and only if $\cP_F(\bv')$ is nonempty.

If $X=X_m$, and we consider $L$-prioritary sheaves, then the above discussion remains true. In particular, if $\bv$ is the Chern character of a type $2$ elementary transformation of $\cV = \cO_{X_m}(-2L)^{\oplus a} \oplus \cO_{X_m}(-L)^{\oplus b}$ along $\bigoplus_{i=1}^{m}\cO_{E_i}(-1)^{\oplus d_{i}}$, then there is an open dense substack of $\cR_L^E(\bv)$ that is an iterated Grassmannian over the point corresponding to $\cO_{\PP^2}^{\oplus a} \oplus \cO_{\PP^2}(L)^{\oplus b}$.
\end{rmk}

We compare the discriminants of the sheaves appearing in Proposition \ref{elementarytransformationexact}.

\begin{prop}
\label{discriminant upper bound elementary transformation}
Let $\cV$ and $\cW$ be sheaves of rank $r\geq 2$ on a surface $X$ containing a smooth rational curve $E$ such that $E^2 = -1$. Suppose $\cV|_E \cong \cO_E^{\oplus r}$, and let $d$ be an integer such that $0<d< r$. If $\cW$ belongs to either of the exact sequences
\[
0 \to \cW(nE) \to \cV((n+1)E) \to \cO_{E}(-n-1)^d \to 0
\]
or
\[
0 \to \cV((n-1)E) \to \cW(nE) \to \cO_{E}(-n)^d \to 0
\]
for any integer $n$, then
\[
\Delta(\cW) = \Delta(\cV) + \frac{d(r-d)}{2r^2}.
\]
In particular, $\Delta(\cW) \leq \Delta(\cV) + 1/8$.
\end{prop}
\begin{proof}
The last statement follows from the first by maximizing the last term in the equality.

Without loss of generality, we can assume $n=0$. The proofs for both exact sequences are nearly identical, so we only prove the equality for the first sequence. We get $\ch(\cO_{E}(-1)) = (0,E,-1/2)$ from the restriction sequence
\[
0 \to \cO \to \cO(E) \to \cO_{E}(-1) \to 0.
\]
From the first sequence, we obtain
\begin{align*}
\Delta(\cW) & = \frac{\ch_1(\cW)^2}{2r^2} - \frac{\ch_2(\cW)}{r}\\
& =  \frac{(\ch_1(\cV(E)) - dE)^2}{2r^2} - \frac{(\ch_2(\cV(E))+ d/2)}{r}\\
& = \Delta(\cV(E)) + \frac{2rd-d^2}{2r^2} - \frac{d}{2r}\\
& = \Delta(\cV(E)) + \frac{d(r-d)}{2r^2}\\
& = \Delta(\cV) + \frac{d(r-d)}{2r^2}.
\end{align*}

\end{proof}

\begin{proof}[Proof of Proposition \ref{refined Bogomolov}]
Twists of $L$-prioritary sheaves are $L$-prioritary, so it is enough to prove the statement for $\cV$.

Let $\pi:X_m \to \PP^2$ be the map contracting the exceptional divisors $E_i$. By our assumption, the bundle $\cV$ belongs to an exact sequence
\[
0 \to \cO_{X_m}(-2L)^{\oplus a} \oplus \cO_{X_m}(-L)^{\oplus b} \to \cV \to \bigoplus_{i=1}^m \cO_{E_i}(-1)^{d_i} \to 0,
\]
and this gives an isomorphism
\[
\pi_*(\cV) \cong \cO_{\PP^2}(-2L)^{\oplus a} \oplus \cO_{\PP^2}(-L)^{\oplus b}.
\]

Suppose towards a contradiction that there is an $L$-prioritary sheaf $\cW$, with $r(\cW) = r(\cV)$, $\nu(\cW) = \nu(\cV)$, and $\Delta(\cW) < \Delta(\cV)$. If $\bw = \ch(\cW)$, then $\cP_{X_m,L}(\bw)$ is nonempty. Furthermore, a general sheaf $\cE \in \cP_{X_m,L}(\bw)$ has rigid splitting type on all exceptional divisors $E_i$ by the discussion in Remark \ref{Grassmannian bundle over stack}. By applying Proposition \ref{elementarytransformationexact}, there is an exact sequence
\[
0 \to \pi^*\pi_*(\cE) \to \cE \to \bigoplus_{i=1}^m \cO_{E_i}(-1)^{\oplus d_i} \to 0,
\]
where $\pi_*(\cE)$ is $L$-prioritary.

By the second sequence in Proposition \ref{discriminant upper bound elementary transformation}, we have
\begin{align*}
\Delta(\pi_*(\cE)) + \sum_{i=1}^m\frac{d_i(r-d_i)}{2r^2} & = \Delta(\cE)\\
& < \Delta(\cV)\\
& = \Delta(\pi_*(\cV)) + \sum_{i=1}^m\frac{d_i(r-d_i)}{2r^2}.
\end{align*}
Thus, we obtain the inequality $\Delta(\pi_*(\cE)) < \Delta(\pi_*(\cV))$, with $r(\pi_*(\cE)) = r(\pi_*(\cV))$ and $\nu(\pi_*(\cE)) = \nu(\pi_*(\cV))$. But this contradicts Proposition \ref{No L-prioritary sheaves with smaller discriminant}, so we are done.
\end{proof}

%%%%%%%%%%%%%%%%%%%%%%%%%%%%%%%%%%%%%%%%%%%%%%%%%%%%%%%%%%%%%%%%%%%%%%%%%%%%%

\section{Construction of good direct sums}
\label{good bundle construction and properties}

In this section, our goal is to construct explicit type $2$ elementary transformations of the bundle
\[
\cO_{X_m}(-2L)^{\oplus a} \oplus \cO_{X_m}(-L)^{\oplus b}
\]
along $\bigoplus_{i=1}^m \cO_{E_i}(-1)^{\oplus d_i}$, where $a$, $b$, $d_1$, \ldots, and $d_m$ are integers as in Proposition \ref{refined Bogomolov}, and show their higher cohomology vanishes. A natural approach, which is the one we take, is to ``distribute'' $d_i$ copies of $E_i$ among the summands and show they have positivity properties. We choose $\nu(\cV).L$ to lie in the interval $[-2,-1]$ since it is the minimal total slope with respect to intersecting $L$ that guarantees $H^2(X_m,\cV) = 0$ (Proposition \ref{h0h2vanishing}).

\begin{construction}
\label{bundle construction}
Suppose as above, we start with a bundle $\cV$ and some integers $d_i$. Order the $r$ summands of $\cV$ in an $r$-tuple by increasing $-K_{X_m}$-slope. In this case, it is quite simple, and we have
\[
\cS:=(\underbrace{\cO(-2L),\ldots,\cO(-2L)}_\text{$a$ copies},\underbrace{\cO(-L),\ldots,\cO(-L)}_\text{$b$ copies}).
\]

\begin{enumerate}
        \item[]\textbf{Step 1)} Start with $i=1$. Twist each coordinate by $\cO(E_i)$ starting from left to right in $\cS$ until reaching the $d_i$th coordinate.
        \item[]\textbf{Step 2)} Let $\cS'$ be the new $r$-tuple obtained from the previous step. Reorder the coordinates of $\cS'$ by increasing $-K_X$-slope. If two distinct line bundles $\cO(D_1)$ and $\cO(D_2)$ have the same $-K_X$-slope, then $\cO(D_1)$ sits to the left of $\cO(D_2)$ if either
        \begin{enumerate}
                \item[1)] $-D_1.L < -D_2.L$,
                \item[2)] or $D_1.L=D_2.L$ and there exists a $j$ such that $D_1.E_i=D_2.E_i$ for all $i < j$ and $D_1.E_j < D_2.E_j$.
        \end{enumerate}
        
        \item[]\textbf{Step 3)} Repeat steps 1) and 2) using $E_{i+1}$.
        
\end{enumerate}

We call a bundle $\cV'$ constructed in this way a \textit{good} bundle.

\begin{rmk}
\label{good bundle for every total slope}
Given a rank $r$ and total slope $\nu$, there is a unique (up to isomorphism) good bundle $\cV$ such that $r(\cV) = r$ and $\nu(\cV) = \nu$ by Remark \ref{all total slopes possible}.
\end{rmk}

\end{construction}

\begin{prop}[Lemma 2.2 \cite{WalterComponentsofStack}]
\label{balancedimpliesprior}
Let $X$ be a surface. Let $C \subset X$ be a smooth rational curve such that $h^0(X, \cO(C)) \geq 2$, and let $n = -K_X.C$. If $\cE$ and $\cG$ are torsion-free sheaves whose restrictions to $C$ are $\cE|_C \cong \bigoplus_i \cO_C(e_i)$ and $\cG|_C \cong \bigoplus_j \cO_C(g_j)$ with $\max\{e_i\} - \min\{g_j\} < n$, then $\Ext^2(\cE,\cG) = 0$. In particular, if $\cE$ is a torsion-free sheaf on $X_m$ and $\cE$ is of rigid splitting along some line $\ell \in |L|$, then $\cE$ is $L$-prioritary.
\end{prop}

\begin{proof}
By Serre duality, we have $\Ext^2(\cE,\cG) \cong \Hom(\cG,\cE(K_X))^*$. A morphism $\phi \in \Hom(\cG,\cE(K_X))$ induces a morphism $\phi|_{C} \in \Hom_{\cO_C}(\cG|_C,\cE(K_X)|_C) \cong H^0(C,\bigoplus_{i,j} \cO_C(e_i-g_j-n)) = 0$. Since $C$ is movable, $\phi$ vanishes on all of $S$, and we obtain $\Ext^2(\cE,\cG) = 0$. To prove the last statement, set $C = L$ and $\cG = \cE(-L)$.
\end{proof}

The above shows good bundles are $L$-prioritary. Also, good bundles are type $2$ elementary transformations of $\cV$ by construction.

\begin{exmp}
Here is an example of Construction \ref{bundle construction}. Suppose we want a bundle on $X_3$, and we start with the bundle
\[
\cV = \cO(-2L)^{\oplus 3} \oplus \cO(-L)^{\oplus 2}.
\]
We list the good bundles $\cV'$ for some choices of $(d_1,d_2,d_3)$.

\begin{enumerate}
\item[(1,0,0):] $\cV' = \cO(-2L)^{\oplus 2} \oplus \cO(-2L+E_1) \oplus \cO(-L)^{\oplus 2}$
\item[(1,1,0):] $\cV' = \cO(-2L) \oplus \cO(-2L+E_1) \oplus \cO(-2L+E_2) \oplus \cO(-L)^{\oplus 2}$
\item[(2,1,1):] $\cV' = \cO(-2L+E_1) \oplus \cO(-2L+E_2) \oplus \cO(-2L+E_1+E_3) \oplus \cO(-L)^{\oplus 2}$
\item[(2,2,2):] $\cV' = \cO(-2L+E_1+E_2) \oplus \cO(-2L+E_1+E_3) \oplus \cO(-2L+E_2+E_3) \oplus \cO(-L)^{\oplus 2}$
\item[(3,2,2):] $\cV' = \cO(-2L+E_1+E_2) \oplus \cO(-2L+E_1+E_3) \oplus \cO(-L)^{\oplus 2} \oplus \cO(-2L+E_1+E_2+E_3)$
\item[(3,2,4):] $\cV'=\cO(-2L+E_1+E_3) \oplus \cO(-L) \oplus \cO(-2L+E_1+E_2+E_3)^{\oplus2} \oplus \cO(-L+E_3)$
\end{enumerate}

\end{exmp}

\begin{prop}
\label{summands have small slope difference}
Let $\cV = \bigoplus_{i=1}^r \cO(D_i)$ be a good bundle. Then $|(D_i-D_j).(-K_{X_m})| \leq 3$ with equality if and only if $D_i-D_j = \pm L$. Furthermore, if $D_i.L = D_j.L$, then $|(D_i-D_j).(-K_{X_m})| \leq 1$.
\end{prop}

\begin{proof}
We induct on the number of points of the blowup. For $\PP^2$, the statement is clear. Suppose the statement is true for $X_{k-1}$. Let $\cV' \cong \bigoplus_{i=1}^r \cO(D_i')$ be the pushforward of $\cV$ to $X_{k-1}$.
\begin{enumerate}
\item[1)] Suppose $D_i = D_i'+E_k$ and $D_j=D_j'+E_k$ (or $D_i = D_i'$ and $D_j=D_j'$). Then $D_i-D_j=D_i'-D_j'$, so the statement is true since it is true for $\cV'$ by induction.
\item[2)] Without loss of generality, suppose $D_i = D_i'+E_k$ and $D_j=D_j'$. First, consider the case when $D_i.L = D_j.L$. Then $D_i'.L=D_j'.L$.\\

        \begin{enumerate}
        \item[a)] If $|(D_i'-D_j').(-K_{X_k})|= 0$, then the statement is true.
        \item[b)] Suppose $|(D_i'-D_j').(-K_{X_k})|= 1$. Since $D_i=D_i'+E_k$, the bundle $\cO(D_i')$ precedes $\cO(D_j')$ in the ordering of Step 2) in Construction \ref{bundle construction}, hence has strictly smaller $-K_{X_k}$-slope. Thus, we obtain $(D_i'-D_j').(-K_{X_k}) = -1$ and $(D_i-D_j).(-K_{X_k}) = (D_i'-D_j').(-K_{X_k}) + 1=0$.\\
        \end{enumerate}
        Now assume $D_i.L \neq D_j.L$. By induction, we know $|(D_i'-D_j').(-K_{X_k})| \leq 3$. If $(D_i'-D_j').(-K_{X_k})\neq 0$, then by the same reasoning as in b) above, we have that $\cO(D_i')$ has smaller $(-K_{X_k})$-slope, hence $|(D_i-D_j).(-K_X)| \leq 2$. If $(D_i'-D_j').(-K_X)=0$, then the statement is obvious.
\end{enumerate}
\end{proof}

From now on, we restrict to $X_m$ for $m \leq 5$.

\begin{lemma}
\label{small intersection difference with rigid curves}
Let $\cV = \bigoplus_{i=1}^r \cO(D_i)$ be a good bundle. For any $(-1)$-curve $C$, we have $|(D_i-D_j).C| \leq 2$.
\end{lemma}
Here we really need to restrict to degree at least $4$ del Pezzo surfaces. For example, on $X_6$ we get $|(D_i-D_j).C| =3$ when $C=L-E_4-E_6$ and $\cV = \cO(-2L+E_1+E_2+E_3+E_5) \oplus \cO(-L+E_4+E_6)$.

\begin{proof}
When $C=E_k$, we get the bound $|(D_i-D_j).C| \leq 1$ by construction. When $C$ is not of this form, we use induction on the number of points of the blowup. For $\PP^2$, there is nothing to check. Suppose now it is true for $X_{k-1}$. Let $\pi:X_k \to X_{k-1}$ be the map contracting $E_k$, and let $\pi_*(\cV) =\bigoplus_{i=1}^r \cO(D_i')$. If $D_i-D_j=\pi^*(D_i'-D_j')$, then the bound holds by induction.

To avoid clutter, we will drop the $\pi^*$ notation when referring to $D_i'$ and $D_j'$ as divisors on $X_k$. Without loss of generality, suppose $D_i=D_i'+E_k$ and $D_j=D_j'$. First, consider the case $C=L-E_k-E_l$. By the induction hypothesis, we have $-2 \leq (D_i'-D_j').C \leq 2$. If $-2 \leq (D_i'-D_j').C \leq 1$, then we are done. Thus, we assume towards a contradiction that $(D_i'-D_j').C=2$. Since $(D_i'-D_j').E_k = 0$, we obtain
\begin{align*}
(D_i'-D_j').C & = (D_i'-D_j').(L-E_k-E_l)\\
& = (D_i'-D_j').(L-E_l) - (D_i'-D_j').E_k\\
& = (D_i'-D_j').(L-E_l).
\end{align*}
Note that $(D_i'-D_j').L \in \{-1,0,1\}$ and $(D_i'-D_j').(-E_l) \in \{-1,0,1\}$ by construction of the divisors $D_i'$ and $D_j'$. Thus, the only way to achieve $(D_i'-D_j').(L-E_l) = 2$ is if $(D_i'-D_j').L = 1$ and $(D_i'-D_j').(-E_l) = 1$. From these equalities, we find
\begin{align*}
D_i'.L & = -1,\\
D_i'.E_l & = -1,\\
D_j'.L & = -2,\\
\text{and } D_j'.E_l & = 0.
\end{align*}
We  claim $(D_i'-D_j').(-K_{X_k})=0$. Since $-D_i'.L < -D_j'.L$, the only way for $E_k$ to be added to $D_i'$ before $D_j'$ is in Construction \ref{bundle construction} Step 2), which forces $(D_i'-D_j').(-K_{X_k})=0$. Observe that $D_i'.(-K_{X_k}) \geq -2$ since $D_i'.L = -1$ and $D_i'.E_l = -1$. Therefore, we have $D_j'.(-K_{X_k}) = D_i'.(-K_{X_k}) \geq -2$. However, this is absurd; we must have $D_j'.(-K_{X_k}) \leq -3$ since $D_j'.L = -2$ and $0 \leq D_j'.E_i \leq 1$ for all $i$ with $D_j'.E_k = D_j'.E_l = 0$.

Second, we consider the case $C=2L-\sum_{n=1}^5 E_{n}$. Decompose $-K_{X_5}$ into a sum of rational curves $-K_{X_5} = C + L$. As in the previous case, we must show $(D_i'-D_j').C \neq 2$. If we have equality, then
\[
(D_i'-D_j').(-K_{X_5}) = 2+(D_i'-D_j').L.
\]
Since $D_i=D_i'+E_k$, we must have $(D_i'-D_j').(-K_{X_5}) \leq 0$. This is impossible since $|(D_i'-D_j').(L)| \leq 1$.
\end{proof}

\begin{prop}
\label{special bad case}
Let $\cV = \bigoplus_{i=1}^r \cO(D_i)$ be a good bundle. If $C_1$ and $C_2$ are distinct $(-1)$-curves on $X_m$ and $(D_{i_1}-D_{j}).C_1 = (D_{i_2}-D_{j}).C_2 = 2$ for some $i_1$,$i_2$, and $j$, then $C_1.C_2 = 1$ if and only if $m = 5$, the $(-1)$-curves $C_1$ and $C_2$ are $C_1=L-E_1-E_4$ and $C_2=L-E_2-E_5$, and the divisors are either
\begin{enumerate}
\item[(i)] $D_{i_1} = -2L+E_1+E_4$, $D_{i_2} = -2L+E_2+E_5$, and $D_j=-2L+E_3$,
\item[(ii)] or $D_{i_1} = -L+E_1+E_4$, $D_{i_2} = -L+E_2+E_5$, and $D_j=-L+E_3$.
\end{enumerate}
Furthermore, if $(D_1,D_2,D_3)$ and $(D_4,D_5,D_6)$ are two such triples occurring in $\cV$, then they are identical up to permutation.
\end{prop}

\begin{proof}
Since $(D_i-D_j).E_l \leq 1$ and $(L-E_i-E_j).(2L-\sum_{n=1}^5 E_{n}) = 0$ for any $i$, $j$, and $l$, the only case to check is when $C_1$ and $C_2$ are lines passing through disjoint pairs of exceptional divisors. We eliminate the possibilities case by case. First suppose $D_{i_1}.L=-2$ and $D_{j}.L = -1$. Then $D_{i_1}.C_1 \leq 0$, $D_{j}.C_1 \geq -1$, and $(D_{i_1}-D_j).C_1 \leq 1$, contradicting the assumption on the intersection value. Let $C_i = L-E_{k_i}-E_{l_i}$. If $D_{i_1}.L=-1$ and $D_{j}.L = -2$, then $(D_{i_1}-D_j).C_1 =1-(D_{i_1}-D_j).(E_{k_1}+E_{l_1})$. Without loss of generality, suppose $(D_{i_1}-D_j).E_{k_1}=-1$. Then $D_{i_1} = -L+E_{k_1} + \cdots$. By Construction \ref{bundle construction}, this implies $D_j.(-K_{X_m}) \geq -3$ since $E_k$ is added to a divisor $D$ with $D.L=-1$ before a divisor $D'$ with $D'.L=-2$ if and only if $D'.(-K_{X_m})=D.(-K_{X_m}) \geq -3$. By the same analysis for $C_2$, we have $D_{i_2} = -L+E_{k_2}+ \cdots$, and $E_{k_i}.D_j = 0$ for $i=1,2$. Since $m \leq 5$, we must have $D_j=-2L+E_{l_1}+E_{l_2}+E_n$, where $E_n \neq E_{k_i}$ for both $i$. This now forces $D_{i_s} = -L+E_{k_s}+E_{l_s}+\cdots$, which means $(D_{i_s}-D_j).(-K_{X_m}) \geq 2$. The contradiction arises from the fact that $D_{i_s}.L=-1$, $D_{i_s}.E_{k_s}=-1$, and $D_j.E_{k_s}=0$, but their $-K_{X_m}$-slopes differ by more than $1$. The construction does not allow for this to happen; once two divisors have equal $-K_{X_m}$-slope, they cannot have a difference of slopes strictly greater than $1$. Thus, we must have that $D_{i_1}$,$D_{i_2}$, and $D_j$ have the same intersection with $L$.

We are reduced to showing that the divisors are of the specified type. By Proposition \ref{summands have small slope difference}, we must have $(D_{i_s}-D_j).(-K_{X_m}) \leq 1$. Without loss of generality, suppose $D_j.L=-2$. If $D_j.(-K_{X_m}) \neq -5$, then it is easy to see that the assumed difference of intersection values on the $C_i$ cannot occur. From this, we deduce that $D_j = -2L+E_{n_1}$, $D_{i_1} = -2L+E_{n_2}+E_{n_3}$, and $D_{i_2}=-2L+E_{n_4}+E_{n_5}$. The construction forces the $n_i$ to be the specified values. For the last part of the proposition, note that the construction does not allow for $\cO(-2L+E_1+E_4)$ and $\cO(-L+E_3)$ appearing as summands of $\cV$. 
\end{proof}

\begin{corollary}
\label{no positive contribution to bad case from other divisors}
Suppose $\cV$, $D_{i_1}$, $D_{i_2}$, $D_j$, $C_1$, and $C_2$ are as in Proposition \ref{special bad case}. If $(D_s-D_j).C_l=2$ for some $i$, then $D_s = D_{i_l}$. In particular, if $D_s \neq D_{i_l}$, then $(D_s-D_j).C_l \leq 1$.
\end{corollary}

\begin{proof}
Replace $D_{i_l}$ with $D_s$ and use Proposition \ref{special bad case}. For the last statement, use Lemma \ref{small intersection difference with rigid curves}.
\end{proof}

\begin{prop}
\label{SHGH for good bundles}
Let $\cV = \bigoplus_{i=1}^r \cO(D_i)$ be a good bundle on $X_m$. If $D$ is a divisor such that $\nu(\cV(D))$ is nef, then $H^i(X_m,\cV(D)) = 0$ for $i > 0$.
\end{prop}

\begin{proof}
Since $\nu(\cV(D))$ is nef, we have $(D_j+D).C \geq -1$ for all $(-1)$-curves by Lemma \ref{small intersection difference with rigid curves}. If $(D_j+D).C_1 = (D_j+D).C_2 = -1$, then $C_1.C_2 = 1$ if and only if $C_1 = L + E_1+E_4$, $C_2 = L+E_2+E_5$, and $D_j = -2L+E_3$ or $D_j=-L+E_3$ by Proposition \ref{special bad case}. Without loss of generality, suppose $D_j = -2L+E_3$. By Corollary \ref{no positive contribution to bad case from other divisors}, we have $(D_i+D).C_1 \leq 0$ and $(D_i+D).C_2 \leq 0$ unless $D_i = -2L+E_1+E_4$ or $D_i=-2L+E_2+E_5$. Using the notation from Proposition \ref{special bad case}, we set $D_{i_1} = -2L+E_1+E_4$ and $D_{i_2} = -2L+E_2+E_5$. Note that $D.C_1 = D.C_2 = 1$, so $(D_{i_1}+D).C_1 = 1$, $(D_{i_1}+D).C_2 = -1$, $(D_{i_2}+D).C_1 = -1$, and $(D_{i_2}+D).C_2 = 1$. If $\cO(D_j + D)$, $\cO(D_{i_1}+D)$, and $\cO(D_{i_2}+D)$ appear in $\cV(D)$ as direct summands $a$, $b$, and $c$ times, respectively, then $0 \leq \ch_1(\cV(D)).C_1 \leq b-a-c$ and $0 \leq \ch_1(\cV(D)).C_2 \leq c-a-b$. From these inequalities, we obtain $a+b \leq c$ and $a+c \leq b$. Combining these two inequalities gives $a+c \leq c-a$, which shows $a = 0$. Thus, whenever $\ch_1(\cV(D))$ is nef and $m \geq 2$, all of the summands of $\cV(D)$ satisfy the hypotheses of Corollary \ref{line bundle cohomology}. For $m \leq 1$, it is easy to check directly that the summands have no higher cohomology.
\end{proof}

%%%%%%%%%%%%%%%%%%%%%%%%%%%%%%%%%%%%%%%%%%%%%%%%%%%%%%%%%%%%%%%%%%%%%%%%%%%%%%%%%%%%%

\section{Brill-Noether}
\label{Brill-Noether}
Throughout this section, we work with del Pezzo surfaces of degree at least $4$. Let $W_m$ be the Weyl group that acts on $\Pic(X_m)$. Observe that $K_{X_m}+\sigma(L)$ is anti-effective for all $\sigma \in W_m$, so $H.(K_{X_m}+\sigma(L)) < 0$. If $m \geq 2$ and $C$ is any $(-1)$-curve, then there exists some $\sigma$ such that $\sigma(L)-C$ is a fiber class. Therefore, all $\mu_H$-semistable sheaves are $\sigma(L)$-prioritary and $C$-prioritary for all $\sigma \in W_m$ and $(-1)$-curves $C$ by Proposition \ref{effective+prioritary}. We will only need the following consequence of the existence of $\mu_H$-semistable sheaves in $\cP_L(\bv)$: If there exists a $\mu_H$-semistable sheaf in $\cP_L(\bv)$ and $\cE$ is a general sheaf in $\cP_L(\bv)$, then $\cE$ has rigid splitting type on all $(-1)$-curves by Lemma \ref{prioritary smooth restriction morphism}. For this section, it is enough to work with the stack of sheaves $\cP_{L^W}(\bv)$ that are prioritary with respect to $\sigma(L)$ for all $\sigma \in W_m$.

We compute the cohomology of a general sheaf in $\cP_{L^W}(\bv)$. Since a general sheaf in $\cP_{L^W}(\bv)$ is locally free by \cite{WalterIrreducibility} and the Serre dual of a prioritary bundle is prioritary, it makes no difference to consider $\bv$ or its Serre dual $\bv^* \otimes K_X$. Thus, we can assume $\nu(\bv).L \geq -2$.

Consider the map $\pi:X_m \to X_{m-1}$ contracting $E_m$. For any $\bv \in K(X_m)$ such that $\cP_{L^W}(\bv) \neq \emptyset$, a general sheaf $\cE \in \cP_{L^W}(\bv)$ is of rigid splitting type along $E_m$. If $\cE$ and $\cF$ are sheaves in $\cP_{L^W}(\bv)$ of rigid splitting type along $E_m$ and $\nu(\bv).E_m \leq 0$, then $[R^1\pi_*(\cE)] = [R^1\pi_*(\cF)]$ and  $[\pi_*(\cE)] = [\pi_*(\cF)]$ in $K(X_{m-1})$. We refer to the class of $[\pi_*(\cE)]$ as $\pi_*(\bv)$. Our main theorem is the following.

\begin{theorem}[Brill-Noether]
\label{Brill-Noether theorem}
Let $X_m$ be a del Pezzo surface with $m \leq 5$. Let $\bv$ be a Chern character such that $\nu(\bv).L \geq -2$ and $\cP_{L^W}(\bv) \neq \emptyset$.
\begin{enumerate}
\item[(1)] If $\nu(\bv).C \geq -1$ for all $(-1)$-curves $C$, then $\bv$ is non-special.
\item[(2)] Let $W_m$ be the Weyl group acting on $\Pic(X_m)$. If there exists some $\sigma \in W_m$ such that $\nu(\bv).\sigma(L) \leq -1$ or $\nu(\bv).\sigma(L-E_i) \leq -1$ for some $i$, then $\bv$ is non-special.
\item[(3)] Let $\nu(\bv).\sigma(L) > -1$ and $\nu(\bv).\sigma(L-E_1) > -1$ for all $\sigma \in W_m$. Suppose $C$ is a $(-1)$-curve such that $\nu(\bv).C < -1$ and let $\pi$ be the map contracting $C$. Then $\bv$ is non-special if and only if $\pi_*(\bv)$ is non-special on $X_{m-1}$ (or $\PP^1 \times \PP^1$) and $\chi(\pi_*(\bv)) \leq 0$. 
\end{enumerate}
\end{theorem}

\begin{rmk}
By Grothendieck-Riemann-Roch, we have 
\[
\ch(\pi_!(\bv)) = \left(r(\bv), \hspace{2mm} \pi_*(\ch_1(\bv)), \hspace{2mm} \ch_2(\bv)-\frac{E_m.\ch_1(\bv)}{2}\right).
\]
The sheaf $R^1\pi_*(\cE)$ is supported at a single point, so $\ch(R^1\pi_*(\cE)) = (0,0,h^0(X_{m-1},R^1\pi_*(\cE)))$. Therefore, the Chern character of $\pi_*(\bv)$ is
\[
\ch(\pi_*(\bv)) = \left(r(\bv), \hspace{2mm} \pi_*(\ch_1(\bv)), \hspace{2mm} \ch_2(\bv)-\frac{E_m.\ch_1(\bv)}{2} + h^0(X_{m-1},R^1\pi_*(\cE))\right).
\]
We will show how to compute $h^0(X_{m-1},R^1\pi_*(\cE))$ in Section \ref{nef case Brill-Noether}.
\end{rmk}

When $\bv$ is the Chern character of a good bundle and $\ch_1(\bv)$ is close to being a nef class, we show that the higher cohomology of a general sheaf vanishes, which proves Conjecture \ref{Higher rank SHGH} for del Pezzo surfaces of degree at least $4$. We obtain the non-speciality of $L$-prioritary sheaves with higher discriminants by using elementary modifications discussed below. If $\ch_1(\bv)$ is far from nef, then we contract the $(-1)$-curves $C$ such that $\nu(\bv).C < -1$ and compare the cohomology of a general sheaf $\cE \in \cP_{L^W}(\bv)$ with the cohomology of $\pi_*(\cE)$.

First, we give a simple criteria for the vanishing of $H^2(X_m,\cE)$ for a general sheaf $\cE \in \cP_{L^W}(\bv)$.

\begin{prop}
\label{h2 vanishes for general sheaf}
Suppose $\bv$ is a Chern character on $X_m$ such that $\nu(\bv).L \geq -2$. If $\cE$ is a general sheaf in $\cP_{L}(\bv)$, then $H^2(X_m,\cE) = 0$. In particular, a general sheaf $\cE$ in $\cP_{L^W}(\bv)$ satisfies $H^2(X_m,\cE) = 0$.
\end{prop}

\begin{proof}
For any divisor $D$ satisfying $D.L \geq 0$ and good bundle $\cV$, the summands of $\cV(D)$ have no cohomology in degree $2$ by Proposition \ref{h0h2vanishing}. We conclude by Proposition \ref{refined Bogomolov} and Lemma \ref{elementarymod}.
\end{proof}

We reproduce the statement for the classification of non-special Chern characters for Hirzebruch surfaces since we will be referring to it frequently.

\begin{theorem}[Corollary 3.7, Corollary 3.9 \cite{CoskunHuizengaBrillNoetherGloballyGenerated}]
\label{Brill-Noether for Hirzebruchs}
Let $\FF_e$ be a Hirzebruch surface. Let $\bv \in K(\FF_e)$ be a Chern character of positive rank with $\Delta(\bv) \geq 0$. Let $\cE \in \cP_F(\bv)$ be a general sheaf.
\begin{enumerate}
\item[(1)] If $\nu(\bv).F \geq -1$, then $H^2(\FF_e, \cE) = 0$.
\item[(2)] If $\nu(\bv).F \leq -1$, then $H^0(\FF_e, \cE) = 0$.
\end{enumerate}
Assume $\nu(\bv).F \geq -1$. Then $\bv$ is non-special if and only if one of the following holds.
\begin{enumerate}
\item[(3)] We have $\nu(\bv).F = -1$.
\item[(4)] We have $\nu(\bv).F > -1$ and $\nu(\bv).E \geq -1$.
\item[(5)] If $\nu(\bv).F > -1$ and $\nu(\bv).E < -1$, let $ m$ be the smallest positive integer such that either $\nu(\bv(-mE)).F \leq -1$ or $\nu(\bv(-mE)) \geq -1$.
        \begin{enumerate}
        \item[(a)] If $\nu(\bv(-mE)).F \leq -1$, then $\bv$ is non-special.
        \item[(b)] If $\nu(\bv(-mE)).F >\ -1$, then $\bv$ is non-special if and only if $\chi(\bv(-mE)) \leq 0$.
        \end{enumerate}
        \qed
\end{enumerate}
\end{theorem}

\subsection{Elementary modifications and Brill-Noether}
\label{nef case Brill-Noether}
We will need to introduce elementary modifications in order to prove Theorem \ref{Brill-Noether theorem} (1).

\begin{mydef}
Let $X$ be a surface, and let $\cE$ be a torsion-free sheaf on $X$. An \textit{elementary modification of $\cE$} is a coherent sheaf $\cE'$ that is the kernel of a surjective map $\cE \to \cO_p$ for some skyscraper sheaf $\cO_p$.
\end{mydef}

\begin{lemma} [Lemma 2.7 \cite{CoskunHuizengaBrillNoetherGloballyGenerated}]
\label{elementarymod}
Let $D$ be a divisor on a surface $X$, and consider a general elementary modification $\cE'$ of a torsion-free sheaf $\cE$ at a general point $p$. Then:
\begin{enumerate}
\item[(1)] If $\cE$ is $D$-prioritary, then $\cE'$ is $D$-prioritary.
\item[(2)] $(r(\cE'),\ch_1(\cE'),\ch_2(\cE')) = (r(\cE),\ch_1(\cE),\ch_2(\cE) - 1)$, and 
\begin{gather*}
\chi(\cE') = \chi(\cE)-1,\\
\Delta(\cE') = \Delta(\cE) + \frac{1}{r(\cE)}.
\end{gather*}
\item[(3)] $H^2(X,\cE') \cong H^2(X,\cE)$.
\item[(4)] If at least one of $H^0(X,\cE)$ or $H^1(X,\cE)$ vanishes, then at least one of $H^0(X,\cE')$ or $H^1(X,\cE')$ vanishes. In particular, if $\cE$ is non-special and $H^2(X,\cE) = 0$, then the same is true for $\cE'$.
\end{enumerate}
\end{lemma}

\begin{rmk}
\label{Sheaves of small discriminant produce sheaves of larger discriminant}
Suppose $\bv$ and $\bv'$ are two Chern characters such that $\ch_0(\bv) = \ch_0(\bv') \geq 1$ and $\ch_1(\bv) = \ch_1(\bv')$. If $\Delta(\bv) < \Delta(\bv')$, then $\Delta(\bv') - \Delta(\bv) = m/\ch_0(\bv)$ for some integer $m > 0$. Indeed, if we expand out the formulas for $\Delta$ we find
\[
\Delta(\bv') - \Delta(\bv) = \frac{\ch_2(\bv)-\ch_2(\bv')}{\ch_0(\bv)} = \frac{c_2(\bv') - c_2(\bv)}{\ch_0(\bv)} > 0.
\]
Thus, if we can construct an $L$-prioritary sheaf of Chern character $\bv$, then we can always apply elementary modifications a finite number of times to produce an $L$-prioritary sheaf of Chern character $\bv'$. A simple consequence of this is the following.
\end{rmk}

\begin{lemma}
\label{Brill-Noether for low disc implies for high disc}
Let $X_m$ be a del Pezzo surface. Let $\bv$, and $\bv'$ be Chern characters with the same rank and total slope, but with $\Delta(\bv) < \Delta(\bv')$. If $\bv$ is non-special and $\nu(\bv).L \geq -2$, then $\bv'$ is non-special.
\end{lemma}

\begin{proof}
Let $\cE$ be a non-special sheaf in $\cP_{L^W}(\bv)$. Then by Remark \ref{Sheaves of small discriminant produce sheaves of larger discriminant}, a finite number of general elementary modifications of $\cE$ will produce a sheaf $\cE'$ with $\ch(\cE') = \ch(\bv')$. By Proposition \ref{h2 vanishes for general sheaf} and Lemma \ref{elementarymod}, the sheaf $\cE'$ is non-special and $\sigma(L)$-prioritary for all $\sigma \in W_m$.
\end{proof}

\begin{proof}[Proof of Theorem \ref{Brill-Noether theorem} (1)]
We can find a good bundle $\cV = \bigoplus \cO(D_i)$ together with a divisor $D$ such that $r(\cV) = r(\bv)$ and $\nu(\cV (D)) = \nu(\bv)$ (Remark \ref{good bundle for every total slope}). Let $\bw$ be the Chern character such that $\ch(\cV(D)) = \ch(\bw)$. If $\bw'$ is any Chern character such that 
\begin{align*}
r(\bw') & = r(\bw')\\
\nu(\bw') & = \nu(\bw)\\
\Delta(\bw') & < \Delta(\bw),
\end{align*}
then $\cP_L(\bw') = \emptyset$ by Proposition \ref{refined Bogomolov}. Consequently, the stack $\cP_{L^W}(\bw')$ is empty as well. First assume $\nu(\bv)$ is nef. Then $\cV(D)$ has no higher cohomology by Proposition \ref{SHGH for good bundles}. By Lemma \ref{Brill-Noether for low disc implies for high disc}, any Chern character $\bv'$ satisfying $r(\bv')=r(\bv)$, $\nu(\bv') = \nu(\bv)$, and $\Delta(\bv') \geq \Delta(\bv)$ is non-special. Thus, Theorem \ref{Brill-Noether theorem} (1) is proved when $\nu(\bv)$ is nef.

If $-1 \leq \nu(\bv).C < 0$ for some $C$, then we can consider the map $\pi:X_m \to X_{m-1}$ (or possibly $\pi:X_2 \to \PP^1 \times \PP^1$) contracting $C$. Thus, the higher direct images $R^i\pi_*(\cE)$ vanish for all $i > 0$ by the proof of Proposition \ref{higherdirectimages}. The cohomology of $\cE$ is the cohomology of $\pi_*(\cE)$. We can continue contracting such curves $C$ and taking direct images until we arrive at a sheaf whose first Chern class intersects all $(-1)$-curves nonnegatively. Recall that the rational morphism $\pi_*: \cP_{L^W}(\bv) \dashrightarrow \cP_{L^W}(\pi_*(\bv))$ defined by $\cF \mapsto \pi_*\cF$ is dominant by Remark \ref{Grassmannian bundle over stack}. Now suppose $\psi : X_m \to X_n$ (or $\psi : X_m \to \PP^1 \times \PP^1$) is a map contracting some collection of $(-1)$-curves for which $-1 \leq \nu(\bv).C < 0$ and $\nu(\psi_*(\bv)).C \geq 0$ for the remaining $(-1)$-curves on $X_n$ (a priori there could be many such maps). If $n \geq 2$, then $\nu(\psi_*(\bv))$ is nef, so the general direct image is non-special by the observations above. 

If $n = 1$ and $m \geq 2$, then the map $\psi$ factors through a map $\pi': X_2 \to X_1$. Without loss of generality, we may assume $\psi = \pi'$ and we have $-1 \leq \nu(\bv).E_2 < 0$ and $\nu(\bv).E_1 \geq 0$. If we further assume $-1 \leq \nu(\bv).(L-E_1-E_2) \leq 0$, then we may instead choose $\psi$ to be the map $\psi':X_2 \to \PP^1 \times \PP^1$ contracting $L-E_1-E_2$. Every Chern character arising in this way on $\PP^1 \times \PP^1$ is non-special by Theorem \ref{Brill-Noether for Hirzebruchs}, so we conclude by the observations above. Thus, we may assume $\nu(\bv).(L-E_1-E_2) > 0$. Then $\nu(\bv).(L-E_1) = \nu(\bv).(L-E_1-E_2) + \nu(\bv).E_2 > -1$ and $\nu(\bv).E_1 > 0$, so $\psi_*(\bv)$ is non-special by Theorem \ref{Brill-Noether for Hirzebruchs} (4). If $n = 1$ and $m = 1$, then $\nu(\bv).E_1 \geq 0$, and we have $\nu(\bv).L \geq -2$ by assumption. If $\nu(\bv).F < -1$, then $\bv$ is non-special by Theorem \ref{Brill-Noether for Hirzebruchs} (2) and Propostion \ref{h2 vanishes for general sheaf}. If $\nu(\bv).F \geq -1$, then $\bv$ is non-special by Theorem \ref{Brill-Noether for Hirzebruchs} (4). Finally, if $n = 0$, then every Chern character arising in this way on $\PP^2$ is non-special (\cite{GottscheHirschowitzBrillNoetherP2}), so we conclude again by the observations above.
\end{proof}

\subsection{The Leray spectral sequence} We will prove Theorem \ref{Brill-Noether theorem} (2) and (3). Then we will apply Theorem \ref{Brill-Noether theorem} (1) to determine when $\chi(\pi_*(\bv))$ in Theorem \ref{Brill-Noether theorem} (3) is non-special (Corollary \ref{absolute Brill-Noether} (3)). Let $\pi:X' \to X$ be the blowup of a surface $X$ a point, and let $E \subset X'$ be the exceptional divisor. Given a sheaf $\cE$ on $X'$ such that $\nu(\cE).E \leq 0$, we will want to compute its cohomology using the Leray spectral sequence
\[
E^{p,q}_2 = H^{p}(X,R^q\pi_*(\cE)) \Rightarrow H^{p+q}(X', \cE).
\]
On the $E_2$-page, the spectral sequence takes on the form

\begin{center}
\begin{tikzpicture}
\matrix (m) [matrix of math nodes,
             nodes in empty cells,
             nodes={minimum width=4ex, minimum height=6ex,
                    text depth=1ex,
                    inner sep=0pt, outer sep=0pt,
                    anchor=base},
             column sep=2ex, row sep=2ex]%
{
q \\
\vdots  & \vdots & \vdots & \vdots & \vdots & \ddots \\

1  &   H^0(X, R^1\pi_*(\cE)) & 0 & 0 & 0 & \cdots \\

0 & H^0(X, \pi_*(\cE)) & H^1(X, \pi_*(\cE)) & H^2(X, \pi_*(\cE)) & 0 & \cdots & \\

 & 0 & 1 & 2 & 3 & \cdots & p \\
};
\draw[-stealth] (m-3-2.south east) -- (m-4-4.north west) node[midway,above] {$\delta$};
\draw[thick] (m-1-1.east) -- (m-5-1.east);
\draw[thick] (m-4-1.south) -- (m-4-5.south);
\end{tikzpicture}
\end{center}
with $\delta$ being the only possible nonzero differential. If $\delta = 0$, then the spectral sequence degenerates on the $E_2$-page, and we get the isomorphisms
\begin{align*}
H^0(X',\cE) & \cong H^0(X,\pi_*(\cE)),\\
H^1(X',\cE) & \cong H^{1}(X,\pi_*(\cE)) \oplus H^0(X, R^1\pi_*(\cE)),\\
H^2(X',\cE) & \cong H^2(X,\pi_*(\cE)).
\end{align*}

In general, we cannot expect the spectral sequence to degenerate on the $E_2$-page. For example, consider the line bundle $\cO_{X_1}(-3L+2E_1)$. It is easy to see $\pi_*\cO_{X_1}(-3L +2E_1) = \cO_{\PP^2}(-3L)$ and $R^1\pi_*(\cO_{X_1}(-3L+2E_1))$ is a skyscraper sheaf supported at the center of the blowup, so we have $H^2(\PP^2,\pi_*\cO_{X_1}(-3L +2E_1)) = \CC$ and $H^0(\PP^2,R^1\pi_*(\cO_{X_1}(-3L+2E_1))) = \CC$. Thus, the Leray spectral sequence is not helpful for computing cohomology in this case, and, in fact, we know that the differential $\delta:H^0(\PP^2,R^1\pi_*(\cO_{X_1}(-3L+2E_1))) \to H^2(\PP^2,\pi_*\cO_{X_1}(-3L +2E_1)) = \CC$ is an isomorphism since the cohomology of $\cO_{X_1}(-3L+2E_1)$ can be determined from the cohomology of $\cO_{X_1}(-E_1)$ by Serre duality, and $\cO_{X_1}(-E_1)$ has no cohomology. Fortunately, these cases can be easily handled by Theorem \ref{Brill-Noether theorem} (2).

\begin{proof}[Proof of Theorem \ref{Brill-Noether theorem} (2)]
Suppose there exists some $\sigma \in W_m$ such that $\ch_1(\bv).\sigma(L) \leq -1$. By Proposition \ref{h2 vanishes for general sheaf}, a general sheaf $\cE \in \cP_L(\bv)$ satisfies $H^2(X_m, \cE) = 0$. On the other hand, a general sheaf $\cE \in \cP_{\sigma(L)}(\bv)$ also satisfies $H^0(X_m, \cE) = 0$. Indeed, for the Serre dual sheaf $\cE^*(K_{X_m})$ we have
\begin{align*}
\nu(\cE^*(K_{X_m})).\sigma(L) & = -\nu(\cE).\sigma(L) + K_{X_m}.\sigma(L)\\
& = -\nu(\cE).\sigma(L) -3\\
& \geq -2.
\end{align*}
If $\cE \in \cP_{\sigma(L)}(\bv)$ is general, then $\cE^*(K_{X_m})$ is $\sigma(L)$-prioritary, so $H^2(X_m,\cE^*(K_{X_m})) = 0$ by Proposition \ref{h2 vanishes for general sheaf}. Since $\cP_{L^W}(\bv)$ is nonempty, a general sheaf $\cE$ can only have $H^1(X_m,\cE)$ be nonzero.

Suppose there exists some $\sigma$ such that $\nu(\bv).\sigma(L-E_i) \leq -1$. Then there exists some collection of pairwise disjoint $(-1)$-curves $C_1$, \ldots, and $C_{m-1}$ such that $\sigma(L-E_i).C_i = 0$ for all $i$. For each $C_i$, let $k_i$ be the unique integer such that $-1 < \nu(\bv).C_i + k_i \leq 0$, and let $d_i = r(\bv)(\nu(\bv).C_i+k_i)$. Let  
\[
D = \sum_{i=1}^{m-1}k_iC_i.
\]
If $\pi:X_m \to \PP^1 \times \PP^1$ is the map contracting $C_1,\ldots,C_{m-1}$, and $\cE \in \cP_{L^W}(\bv)$ is general,
then there is an exact sequence
\[
0 \to \pi^*\pi_*(\cE(-D)) \to \cE(-D) \to \bigoplus_{i=1}^{m-1}\cO_{C_i}(-1)^{\oplus d_i} \to 0
\]
by Proposition \ref{elementarytransformationexact}. We have $H^0(\PP^1 \times \PP^1, \pi_*(\cE(-D))) = 0$ by Theorem \ref{Brill-Noether for Hirzebruchs} (1), so $H^0(X_m,\cE(-D)) = 0$. If $k_i \geq 0$, then induction on $a \leq 0$ and the exact sequence
\[
0 \to \cE(-D - (a+1)C_i) \to \cE(-D-aC_i) \to \cO_{C_i}(a-1)^{\oplus d_i} \oplus \cO_{C_i}(a)^{\oplus r-d_i} \to 0
\]
gives $H^0(X_m,\cE(-D+k_iC_i)) = 0$. If $k_i \leq -1$, then use the same exact sequence and induction on $a \geq 0$ to get $H^0(X_m,\cE(-D+k_iC_i)) = 0$. Combining these together gives $H^0(X_m,\cE) = 0$.
\end{proof}

When we are not in the situation of Theorem \ref{Brill-Noether theorem} (2), the differential $\delta$ will vanish.

\begin{lemma}
\label{h2 for direct image vanishes}
Let $\pi: X' \to X$ be the blowup of a smooth projective surface $X$ at a point, and let $E$ be the exceptional divisor. Suppose $\cE$ is a torsion-free coherent sheaf on $X'$ such that
\[
\cE|_{E} \cong \cO_E(-k-1)^{\oplus d} \oplus \cO_E(-k)^{\oplus (r-d)}
\]
for some integer $k > 0$. If $H^2(X',\cE(-kE)) = 0$, then $H^2(X,\pi_*(\cE)) = 0$. 
\end{lemma}

\begin{proof}
Consider the restriction sequence
\[
0 \to \cE(-kE) \to \cE(-(k-1)E) \to \cO_E(-2)^{\oplus d} \oplus \cO_E(-1)^{\oplus (r-d)} \to 0.
\]
Then the application of $\pi_*$ to this exact sequence and induction on $k$ shows $\pi_*(\cE(-kE)) \cong \pi_*(\cE)$. Furthermore, we have $R^1\pi_*(\cE(-kE)) = 0$ by the proof of Proposition \ref{higherdirectimages}. The Leray spectral sequence shows 
\[
H^2(X',\cE(-kE)) = H^2(X,\pi_*\cE) = 0.
\]
\end{proof}

Now we can prove Theorem \ref{Brill-Noether theorem} (3).

\begin{proof}[Proof of Theorem \ref{Brill-Noether theorem} (3)]
If either $m \neq 2$ or $m=2$ and $C \neq L-E_1-E_2$, we can use the Weyl group action and assume $C = E_m$. Let $\cE$ be a general sheaf in $\cP_{L^W}(\bv)$, and let the map $\pi: X_m \to X_{m-1}$ be the blow down of $E_m$ (or $\pi: X_2 \to \PP^1 \times \PP^1$ if $m=2$ and $C = L-E_1-E_2$).  Since $E_m.L = 0$, we have $\nu(\cE(-kE_m)).L = \nu(\cE).L \geq -2$ so $H^2(X_m,\cE(-kE_m)) = 0$ by Proposition \ref{h2 vanishes for general sheaf}. Thus, the Leray spectral sequence
\[
E^{p,q}_2 = H^p(X_{m-1}, R^q\pi_*(\cE)) \Rightarrow H^{p+q}(X_m,\cE)
\]
degenerates at the $E_2$-page by Lemma \ref{h2 for direct image vanishes}. This induces isomorphisms
\begin{align*}
H^0(X_m,\cE) & \cong H^0(X_{m-1},\pi_*(\cE))\\
H^1(X_m,\cE) & \cong H^{1}(X_{m-1},\pi_*(\cE)) \oplus H^0(X_{m-1}, R^1\pi_*(\cE))\\
H^2(X_m,\cE) & = 0.
\end{align*}
We see that $\cE$ is non-special if and only if $\pi_*\cE$ is non-special and $\chi(\pi_*(\cE)) \leq 0$.
\end{proof}

\begin{rmk}
The proof of Theorem \ref{Brill-Noether theorem} (3) gives an inductive procedure for computing the cohomology of a general sheaf in $\cP_{L^W}(\bv)$. We can give precise numerical criteria for when $\bv$ is non-special in terms of $\ch(\bv)$.

Let $\pi: X \to Y$ be the blowup of a surface $Y$ at a point $p$ with exceptional divisor $C$. Let $\cE$ be a sheaf on $X$ of rigid splitting type along $C$. Following the proof in Proposition \ref{higherdirectimages}, we see that
\[
R^1\pi_*(\cE) = R^1\pi_*(\cE)_p^{\wedge} = \varprojlim_n H^1\left(C_{n},\cE|_{C_n}\right),
\]
and there are exact sequences
\[
0 \to \cE|_{{C_n}}(n) \to \cE|_{{C_{n+1}}} \to \cE|_{{C_n}} \to 0
\]
for each $n > 0$. Let $k$ be the unique integer such that $-1 \leq (\nu(\cE) - kC).C < 0$ and let $d = (\ch_1(\cE) - r(\cE)kC).C$. Since $\cE$ is general, it has rigid splitting type along $C$, so when $n \geq k$ we have $H^1({C_n},\cE|_{{C_n}}(n))=0$ and the dimension of $H^1({C_n},\cE|_{C_n})$ stabilizes. Furthermore, the maps 
\[
H^1({C_{n+1}},\cE|_{C_{n+1}}) \to H^1({C_n},\cE|_{C_n})
\]
always surject, so $R^1\pi_*(\cE)_p = \CC^m_p$, where $m = \dim_{\CC}(H^1(C_{k-1},\cE|_{C_{k-1}}))$, which can be computed to be
\[
m = \frac{dk(k-1)}{2} + \frac{k(r-d)(k+1)}{2}.
\]
If $H^2(Y, \pi_*\cE) = 0$, then the spectral sequence
\[
E^{p,q}_2 = H^p(Y, R^q\pi_*(\cE)) \Rightarrow H^{p+q}(X,\cE),
\]
degenerates on the $E_2$-page so $\chi(\pi_*(\bv)) = \chi(\bv)+m$. Hence, we have $\chi(\pi_*(\bv)) \leq 0$ if and only if $\chi(\bv) \leq -m$.
\end{rmk}

Let us fix some notation. For any $(-1)$-curve $C$, let $k_C$ be the unique integer such that $-1 \leq (\nu(\bv) - k_CC).C < 0$, and set $-d_C = r(\nu(\bv) - k_CC).C$. When $k_C > 0$, set
\[
m_C = \frac{d_ck_C(k_C-1)}{2} + \frac{k_C(r-d_C)(k_C+1)}{2}.
\]

\begin{corollary}
\label{absolute Brill-Noether}
Let $X$ be a del Pezzo surface of degree at least $4$, and suppose $\bv$ is a Chern character of rank at least $2$ satisfying $\nu(\bv).L \geq -2$ and $\cP_{L^W}(\bv) \neq \emptyset$.
\begin{enumerate}
\item[(1)] If $\nu(\bv).C \geq -1$ for all $(-1)$-curves $C$, then $\bv$ is non-special.
\item[(2)] Let $W_m$ be the Weyl group acting on $\Pic(X_m)$. If there exists some $\sigma \in W_m$ and an integer $1 \leq i \leq m$ such that $\nu(\bv).\sigma(L) \leq -1$ or $\nu(\bv).\sigma(L-E_i) \leq -1$, then $\bv$ is non-special.
\item[(3)] Suppose $\nu(\bv).\sigma(L) \geq -1$ and $\nu(\bv).\sigma(L-E_i) \geq -1$ for all $\sigma \in W_m$ and $i$. Let $C_1$, \ldots, $C_k$ be a collection of disjoint $(-1)$-curves such that
\begin{enumerate}
\item[(i)] $\nu(\bv).C_i < -1$ for all $i$,
\item[(ii)] $C_1$, \ldots, $C_k$ cannot be extended to a longer list of disjoint $(-1)$-curves satisfying (i), and
\item[(iii)] $k$ is minimal among such collections satisfying (i) and (ii).
\end{enumerate}
Then $\bv$ is non-special if and only if $\chi(\bv) \leq -\sum_{i=1}^k m_{C_k}$. 
\end{enumerate}
\qed
\end{corollary}

For purposes of classifying globally generated Chern characters (see \cite{CoskunHuizengaBrillNoetherGloballyGenerated}), it is useful to classify non-special Chern characters with nonnegative Euler characteristic.

\begin{corollary}
\label{positive Euler Brill-Noether X_2}
Let $X_m$ be a del Pezzo surface of degree at least $4$. If $\bv$ is a Chern character of rank at least $2$ with $\cP_{L^W}(\bv) \neq \emptyset$ and $\chi(\bv) \geq 0$, then $\bv$ is non-special if and only if either
\begin{enumerate}
\item[(1)] $\nu(\bv).C \geq -1$ for all $(-1)$-curves $C$,
\item[(2)] Let $W_m$ be the Weyl group acting on $\Pic(X_m)$. There exists a $\sigma \in W_m$ and an  $1 \leq i \leq m$ such that $\nu(\bv).\sigma(L) \leq -1$ or $\nu(\bv).\sigma(L-E_i) \leq -1$.
\end{enumerate}
In particular, if $\chi(\bv) > 0$, then only case (1) is possible.
\end{corollary}

\begin{proof}
In all other cases of Theorem \ref{Brill-Noether theorem}, a general sheaf $\cE$ has nontrivial $H^1(X_m,\cE)$. In (2), the Euler characteristic can only be nonpositive.
\end{proof}

\section{The existence theorem}

Let $(X_m,-K_{X_m})$ be an anti-canonically polarized del Pezzo surface with $m \leq 6$. In this section, we classify the Chern characters $\bv$ such that the moduli space of $-K_{X_m}$-semistable sheaves $M(\bv)$ is nonempty. Theorem \ref{Rudakov's semistable theorem} gives a classification conditional on the existence of \emph{smooth restricted complete families} of sheaves (Definition \ref{restriced and smooth}). We show there exists such families when the Chern character satisfies the DL condition (Definition \ref{DLCondition}), making Theorem \ref{Rudakov's semistable theorem} unconditional.
\subsection{The Dr\'ezet-Le Potier condition}
\label{DL condition and Rudakov theorems}
We recall the main definitions and statements from \cite{RudakovExistenceDelPezzo}.

For a smooth projective variety $X$, an object $\cE \in \cD^b(X)$, the bounded derived category of coherent sheaves on $X$, is called \textit{exceptional} if
\[
\Ext^i_{\cD^b(X)}(\cE,\cE) = \begin{cases}\CC, & i=0\\ 0, &i\neq 0 \end{cases}.
\]
When $X$ is a del Pezzo surface, it is well known that every exceptional object is quasi-isomorphic to a sheaf sitting in some degree. In particular, every torsion-free exceptional sheaf is locally free. For example, line bundles are exceptional objects. See \cite{KuleshovOrlovExceptionals} for a thorough study of exceptional objects on a del Pezzo surface.\\

For the remainder of the section, stability on a del Pezzo surface $X_m$ is always with respect to the anti-canonical polarization $-K_{X_m}$.

\begin{mydef}
\label{DLCondition}
A torsion-free coherent sheaf $\cF$ (or Chern character) satisfies the \textit{Dr\'ezet- Le Potier condition} (abbr. as \textit{DL condition}) if 

\begin{enumerate}[label=DL\arabic*]
\item \label{DL1} for every exceptional bundle $\cE$ satisfying $r(\cE) < r(\cF)$ and
\[
\mu(\cF) \leq \mu(\cE) \leq \mu(\cF) + K_{X_m}^2,
\]
we have $\chi(\cE,\cF) \leq 0$, \\

\item \label{DL2} and  for every exceptional bundle $\cE$ satisfying $r(\cE) < r(\cF)$ and
\[
\mu(\cF)-K_{X_m}^2 \leq \mu(\cE) \leq \mu(\cF),
\]
we have $\chi(\cF,\cE) \leq 0$.
\end{enumerate}

\end{mydef}

Note that \ref{DL1} and \ref{DL2} are equivalent by Serre duality. It is known that the DL condition implies $\Delta(\bv) \geq \frac{1}{2}$ (see \cite{Hirzebruchexist}), so a sheaf $\cF$ satisfying the DL condition is not semi-exceptional.

\begin{prop}
\label{necessaryconditionforstable}
Suppose $\cF$ is a nonexceptional stable sheaf with $r(\cF) \geq 2$, i.e. $\Ext^1(\cF,\cF) \neq 0$. Then $\Delta(\cF) \geq 1/2$, and $\cF$ satisfies the DL condition.
\end{prop}

\begin{proof}
We have $\ext^2(\cF,\cF) = 0$ by Serre duality and stability so $\chi(\cF,\cF) = \hom(\cF,\cF) - \ext^1(\cF,\cF)$. Since $\cF$ is stable, we have $\hom(\cF,\cF) = 1$, and we also have $\ext^1(\cF,\cF) \geq 1$, so
\[
\chi(\cF,\cF) \leq 0.
\]
The inequality follows from the formula $\chi(\cF,\cF) = r(\cF)^2(1-2\Delta(\cF))$. If $\cE$ is an exceptional bundle, then it is stable. Suppose $\cE$ satisfies
\[
\mu(\cF) \leq \mu(\cE) \leq \mu(\cF) + K_{X_m}^2.
\]
Then $\hom(\cE,\cF) = 0$ and $\ext^2(\cE,\cF) = \hom(\cF,\cE(K_{X_m})) = 0$ by stability, which shows $\cF$ satisfies (1) of the DL condition. The other condition is proved similarly.
\end{proof}

\begin{mydef}
\label{restriced and smooth}
Let $\cE$ be a torsion-free sheaf. We say $\cE$ is \textit{restricted} if
\[
\mu_{\max}(\cE) - \mu_{\min}(\cE) \leq K_{X_m}^2,
\]
and $\nu_{\max}(\cE)-\nu_{\min}(\cE) \neq -K_{X_m}$ (Section \ref{stability}).

Let $\cF$ be a flat family of coherent sheaves parametrized by a base $S$.
\begin{enumerate}
\item[(i)] If $S$ is smooth, then $\cF$ is \textit{smooth}.
\item[(ii)] If $\cF_s$ is restricted for every $s \in S$, then $\cF$ is \textit{restricted}.
\end{enumerate}
\end{mydef}

\begin{rmk}
Our definition of a restricted sheaf is slightly more general than the definition in \cite{RudakovExistenceDelPezzo}. The condition that the family of sheaves be restricted in Theorem \ref{Rudakov's semistable theorem} is only needed in Lemma 5.8 \cite{RudakovExistenceDelPezzo}, which can easily be extended to restricted sheaves in our sense.
\end{rmk}

We can now state the converse to Proposition \ref{necessaryconditionforstable}.

\begin{theorem}[4.1, 4.2 \cite{RudakovExistenceDelPezzo}]
\label{Rudakov's semistable theorem}
Let $X_m$ be a del Pezzo surface with anticanonical polarization $-K_{X_m}$, that is not $\PP^1 \times \PP^1$. If $\bv$ is a Chern character satisfying the DL condition and there exists a restricted smooth complete family $\cF$ of sheaves of Chern character $\bv$, then the moduli space $M(\bv)$ is nonempty. In addition, the locus of stable sheaves $M^s(\bv)$ is nonempty if and only if either
\begin{enumerate}
\item[(1)] $\Delta(\bv) = \frac{1}{2}$ and $\bv$ is primitive, or
\item[(2)] $\Delta(\bv) > \frac{1}{2}$.
\end{enumerate}
\end{theorem}

\begin{rmk}
Theorem \ref{Rudakov's semistable theorem} (1) is not stated in \cite{RudakovExistenceDelPezzo}, but can be proved using a standard Shatz strata argument. In the case of $X_1$, the existence of restricted smooth complete families for any Chern character $\bv$ satisfying the DL condition was proved by Rudakov.
\end{rmk}

\begin{theorem}[8.1 \cite{RudakovExistenceDelPezzo}]
\label{Rudakov existence X1}
On $X_1$, Theorem \ref{Rudakov's semistable theorem} holds without the assumption on the existence of restricted smooth complete families. \qed
\end{theorem}

\subsection{Restricted sheaves, smooth complete families, and existence}
\label{restricted sheaves and existence}

Let $X_m$ be a del Pezzo surface of degree at least $3$, and suppose $\bv$ is a Chern character satisfying the DL condition. In this section, we show that there exists a restricted smooth complete family of sheaves of Chern character $\bv$. Our first objective is to show that elementary modifications of restricted sheaves are restricted.\\

The \emph{HN filtration} refers to the Gieseker Harder-Narasimhan filtration, and the \emph{$\mu$-HN filtration} refers to the slope Harder-Narasimhan filtration. We omit the proof of the following standard fact.

\begin{prop}
\label{HNrefinement}
Let $\cE$ be a torsion-free coherent sheaf on a polarized smooth projective variety $(X,H)$. The HN filtration of $\cE$ is a refinement of the $\mu$-HN filtration of $\cE$, i.e. if
\[
0 = \cE_0 \subset \cE_1 \subset \cdots \subset \cE_n = \cE
\]
is the HN filtration of $\cE$ and
\[
0 = \cE_0' \subset \cE_1' \subset \cdots \subset \cE_m' = \cE
\]
is the $\mu$-HN filtration of $\cE$, then there are integers $1 \leq j_1 < j_2 < \ldots < j_n \leq m'$ such that $\cE_{j_i}' = \cE_i$ for all $1 \leq i \leq n$.
\end{prop}

Proposition \ref{HNrefinement} is useful for comparing slopes of the factors in the HN filtration between sheaves and their elementary modifications.

\begin{lemma}
\label{elementarymodificationspreservemaxminslopes}
Let $\cE$ be a torsion-free coherent sheaf on a polarized smooth projective variety $(X,H)$ of dimension at least $2$. Consider the exact sequence
\[
0 \to \cF \to \cE \to \cO_p \to 0
\]
for a general point $p \in X$. Then $\mu_{\max}(\cE) = \mu_{\max}(\cF)$ and $\mu_{\min}(\cE) = \mu_{\min}(\cF)$.
\end{lemma}

\begin{proof}

By Proposition \ref{HNrefinement}, it is enough to work with the slope filtration and slope factors.

Let $\cE_i'$ and $\cF_i'$ be the subsheaves in the $\mu$-HN filtrations of $\cE$ and $\cF$, respectively. We claim that the $\mu$-HN filtrations are related in the following way: 
\begin{enumerate}
\item[(i)] They are of the same length.
\item[(ii)] There exists an integer $k$ such that $\cF_i' = \cE_i'$ for all $i \leq k$.
\item[(iii)] For all $i > k$, the sheaf $\cF_i'$ is an elementary modification of $\cE_i'$.
\end{enumerate}
The statement immediately follows from the claim and the fact that $\mu(\cF) = \mu(\cE)$.

Parts (ii) and (iii) imply (i). Indeed, if the filtration of $\cF$ is of length $a$, then $\cF_a' = \cF$ is an elementary modification of $\cE_a'$. Since $r(\cE_a') = r(\cF) = r(\cE)$, we must have $r(\cE_{a+1}'/\cE_a') = 0$, which means $\cE_a' = \cE_{a+1}' = \cE$. Similarly, by rank considerations $\cE_{a-1}'$ is a proper subsheaf of $\cE$.

We first prove part (ii). Consider the diagram
\[
\begin{tikzcd}
0 \arrow[r] & \cF' \arrow[r] \arrow[d, hook] & \cE_1' \arrow[r] \arrow[d, hook] & \cO_p \arrow[r] \arrow[d,equal] & 0\\
0 \arrow[r] & \cF \arrow[r] & \cE \arrow[r] & \cO_p \arrow[r] & 0.
\end{tikzcd}
\]
Note that the top row is only a complex in general. Since $\mu(\cF_1') = \mu(\cE_1') = \mu(\cF')$, we must have $\cF' \subset \cF_1'$. If $\cE_1' \to \cO_p$ is the zero map, then $\cF' = \cE_1'$ so $\cF_1' \subset \cF'$. If $\cE_1' \to \cO_p$ is surjective, then $\cF_1'$ is a proper subsheaf of $\cE_1'$ so $\cF_1'/\cF'$ is a proper subsheaf of $\cE_1'/\cF' = \cO_p$, but the only such subsheaf is $0$, so we are done.

Now let $k$ be the maximal number such that $\cE_k' \to \cO_p$ is the zero map. By taking quotients, we see $\cF_i' = \cE_i'$ for all $i \leq k$. Furthermore, the quotients $\cF/\cE_k'$ and $\cE/\cE_k'$ combined with the argument above show that $\cF_{k+1}'$ is an elementary modification of $\cE_i$. By taking the cokernels of the vertical arrows in the above diagram, we arrive at the identification $\cF/\cF_{k+1}' \cong \cE/\cE_{k+1}'$, and the $\mu$-HN filtrations are isomorphic. Hence, the subsheaf $\cF_{i}'$ is an elementary modification of $\cE_i'$ for all $i > k$. This proves (iii) and concludes the proof of the lemma.
\end{proof}

\begin{prop}
\label{elementarymodificationsofrestrictedprioritarysheaves}
Let $\cE$ be a restricted $L$-prioritary sheaf. If $\cE'$ is an elementary modification of $\cE$, then $\cE'$ is a restricted $L$-prioritary sheaf. Furthermore, the groups $\Ext^2(\cE,\cE)$ and $\Ext^2(\cE',\cE')$ vanish.
\end{prop}

\begin{proof}
The last statement is a consequence of being $L$-prioritary. The elementary modification $\cE'$ is $L$-prioritary by Lemma \ref{elementarymod} (1), and it is restricted by Lemma \ref{elementarymodificationspreservemaxminslopes}.
\end{proof}

\begin{prop}
\label{constructionofsmoothversalfamilies}
If there exists a restricted $L$-prioritary sheaf of Chern character $\bv$, then there exists a restricted smooth complete family of sheaves of Chern character $\bv$.
\end{prop}

\begin{proof}
Let $\cE$ be a restricted $L$-prioritary sheaf of Chern character $\bv$. Fix an ample $H$, an integer $m \gg 0$, and an integer $r > 0$ such that $\cE$ is quotient of $\cO(-m)^{\oplus r}$. Then there exists a smooth open neighborhood $U \subset \Quot(\cO(-m)^{\oplus r},\bv)$ of $[\cO(-m)^{\oplus r} \to \cE]$ of torsion-free sheaves that is also complete. Indeed, given the exact sequence
\[
0 \to \cK \to \cO(-m)^{\oplus r} \to \cE \to 0
\]
we can apply the functor $\Hom(-,\cE)$. Recall that $\Quot(\cO(-m)^{\oplus r},\bv)$ has a tangent-obstruction theory at $[\cO(-m)^{\oplus r} \to \cE]$ given by the spaces $\Hom(\cK,\cE)$ and $\Ext^1(\cK,\cE)$. Since $m \gg 0$, the groups $\Ext^i(\cO(-m)^{\oplus r},\cE)$ vanish for all $i > 0$. Hence, the Kodaira-Spencer map 
\[
\Hom(\cK,\cE) \to \Ext^1(\cE,\cE)
\]
surjects. Since $\cE$ is $L$-prioritary, we have $\Ext^2(\cE,\cE) = 0$ so $\Ext^1(\cK,\cE)=0$ and $\Quot(\cO(-m)^{\oplus r},\bv)$ is smooth in a neighborhood of $[\cO(-m)^{\oplus r} \to \cE]$. By upper-semicontinuity on the dimension of $\Ext^2(\cE,\cE)$ and openness of being torsion-free, there is a smooth neighborhood $U$ of $[\cO(-m)^{\oplus r} \to \cE]$ that is complete and only contains torsion-free sheaves.

It remains to show that there is a nonempty open subscheme $U' \subset U$ consisting of restricted sheaves. A simple variation of Proposition 2.3.1 \cite{HuybrechtsLehnModuliOfSheaves} shows that restrictedness is an open condition in families, so the set $V \subset U$ consisting of sheaves $\cF$ with $\mu_{\max}(\cF)> \mu_{\max}(\cE)$ or $\mu_{\min}(\cF) < \mu_{\min}(\cE)$ (possibly both) is closed. Thus, we can take $U' = U\setminus V$, which is nonempty since it contains $\cE$.
\end{proof}

\begin{corollary}
\label{small discimrinant rsf implies large discriminant rsf}
Let $\bv$ and $\bv'$ be two Chern characters such that $r(\bv) = r(\bv')$, $\ch_1(\bv) = \ch_1(\bv')$ and $\Delta(\bv) < \Delta(\bv')$. If there is restricted smooth complete family of sheaves of Chern character $\bv$, then there is such a family for $\bv'$.
\end{corollary}

\begin{proof}
By Proposition \ref{constructionofsmoothversalfamilies}, it is enough to find one restricted $L$-prioritary sheaf. We conclude by Proposition \ref{elementarymodificationsofrestrictedprioritarysheaves}. \end{proof}

Thus, we have reduced our problem to executing the following strategy: for every rank and total slope, find a restricted $L$-prioritary sheaf with such numerical data and small discriminant.

The construction of slope filtrations commutes with twists by line bundles so the difference $\mu_{\max}-\mu_{\min}$ is invariant under twists by line bundles. By Lemma \ref{elementarymodificationspreservemaxminslopes}, it is enough to consider slope filtrations. Thus, to show that a restricted sheaf of some fixed rank and total slope exists, it is enough to fix a region in $\Pic(X_m)_{\QQ}$ that tiles all of $\Pic(X_m)_{\QQ}$ under integer translations and solve the problem for the translated slope. By Proposition \ref{summands have small slope difference}, this is achieved by Construction \ref{bundle construction}.

\begin{lemma}
\label{failure of DL condition smaller discriminant}
Let $\bv$ be a Chern character satisfying the DL condition. Let $\cV$ be a coherent sheaf such that $\ch_0(\cV) = \ch_0(\bv)$ and $\nu(\cV)=\nu(\bv)$. If $\cE$ is an exceptional bundle as in (1) (or (2)) of the DL condition and $\chi(\cE,\cV) \geq 0$ (or $\chi(\cV,\cE) \geq 0$), then $\Delta(\cV) \leq \Delta(\bv)$.
\end{lemma}

\begin{proof}
The proof is a simple calculation. We have
\[
\chi(\cE,\cV) \geq 0 \geq \chi(\cE,\bv),
\]
so
\begin{align*}
\chi(\cE,\cV) & = r(\cE)r(\cV)(P(\nu(\cV) - \nu(\cE))-\Delta(\cE)-\Delta(\cV))\\
& \geq r(\cE)r(\bv)(P(\nu(\bv) - \nu(\cE))-\Delta(\cE)-\Delta(\bv))\\
& = \chi(\cE,\bv).
\end{align*}
By the assumptions of $\cV$, the inequality simplifies to
\[
-\Delta(\cV) \geq -\Delta(\bv).
\]
The calculation for (2) is similar.
\end{proof}

\begin{theorem}
\label{existence of rsf}
Let $X_m$ be a del Pezzo surface with $m \leq 6$. If $\bv$ is a Chern character satisfying the DL condition, then there exists a restricted smooth complete family of torsion-free coherent sheaves with Chern character $\bv$.
\end{theorem}

\begin{proof}
Good bundles are $L$-prioritary, and they are restricted by \ref{summands have small slope difference}. Every rank and total slope is achieved by a twist of a good bundle. Thus, by Proposition \ref{constructionofsmoothversalfamilies} and Corollary \ref{small discimrinant rsf implies large discriminant rsf}, there are restricted smooth complete families for all Chern characters with discriminant greater than the discriminant of a good bundle of the same total slope and rank. 

Let $\cV$ be a good bundle. By Lemma \ref{failure of DL condition smaller discriminant}, it is enough to find an exceptional bundle $\cE$ satisfying \ref{DL1} (or \ref{DL2} resp.) such that $\chi(\cE, \cV) \geq 0$ (or $\chi(\cV,\cE) \geq 0$ resp.). If $m \leq 3$, then we can choose $\cE = \cO(K_{X_m})$ for \ref{DL2}. Indeed, we have $-6 \leq \mu(\cV) \leq 0$ since each summand has slope lying in this range, so
\[
\mu(\cV)-K_{X_m}^2 \leq -K_{X_m}^2 \leq \mu(\cV),
\]
and we also have $\chi(\cV,\cO(K_{X_m}))=0$.

If $m \geq 4$, we can use $\cE = \cO(K_{X_m})$ if $m-9 \leq \mu(\cV) \leq 0$. However, we also have to consider the cases 
\begin{enumerate}
\item[(i)] $-6 \leq \mu(\cV)< m-9$,
\item[(ii)] and $0 < \mu(\cV)\leq m-3$.
\end{enumerate}
In (i), we can use $\cE=\cO(-L)$ and \ref{DL1}. For (ii) we can use $\cE=\cO(-2L+\sum_{i=1}^m E_i)$ and \ref{DL2}.
\end{proof}

As a consequence, we obtain:

\begin{theorem}
\label{existence del Pezzo}
Let $X_m$ be a del Pezzo surface of degree at least $3$ with anti-canonical polarization $-K_{X_m}$ that is not $\PP^1 \times \PP^1$. If $\bv$ is a Chern character satisfying the DL condition, then the moduli space $M(\bv)$ is nonempty. In addition, the locus of stable sheaves $M^s(\bv)$ is nonempty if and only if either
\begin{enumerate}
\item[(1)] $\Delta(\bv) = \frac{1}{2}$ and $\bv$ is primitive, or
\item[(2)] $\Delta(\bv) > \frac{1}{2}$.
\qed
\end{enumerate}
\end{theorem}

\begin{rmk} If $\Delta(\bv) \geq 1$, then $\bv$ always satisfies the DL condition so $M^s(\bv)$ is nonempty. Indeed, if $\bv$ and $\cE$ are as in (1) of the DL condition, we have
\[
P(\nu(\bv)-\nu(\cE)) \leq \Delta(\cE)+\Delta(\bv) \iff \chi(\cE,\bv) \leq 0.
\]
Let $\nu(\bv)-\nu(\cE) = a(-K_{X_m})+bD$, where $D.(-K_{X_m}) = 0$. We have $-1\leq a \leq 0$ since $\cE$ is an exceptional bundle in (1) of the DL condition.\ By the Hodge index theorem, we obtain $D^2 \leq 0$, so $(\nu(\bv)-\nu(\cE))^2+(\nu(\bv)-\nu(\cE)).(K_{X_m}) = (a^2+a)(K_{X_m})^2 + b^2D^2 \leq 0$. Thus, we obtain $P(\nu(\bv)-\nu(\cE)) \leq 1$. The argument for (2) is proved similarly.
\end{rmk}

\begin{exmp}
\label{bad semistable example}
There are examples of moduli spaces containing only strictly semistable sheaves that are not included in the theorem. For example, on $X_1$ one can take
\[
\cV_k = \cO(L+E_1)^{\oplus k} \oplus \cO(2L-2E_1).
\]
It is easy to check that $\cV_k$ is Gieseker semistable and $\Delta(\cV_k) \to 0$ as $k \to \infty$. When $k = 1$, we have $\Delta(\cV_1) =1$, so the DL condition is automatically satisfied. For $k \geq 6$, $\Delta(\cV_k) < \frac{1}{2}$, so the moduli space is purely semistable. Let $\bv_k=\ch(\cV_k)$. Then $\dim(M(\bv_k)) \geq 1$ since the stable objects in $M(\bv_1)$ (which form a $4$-dimensional family) produce different $S$-equivalence classes in $M(\bv_k)$.
\end{exmp}

%\bibliography{Brill-Noether2}
%\bibliographystyle{alpha}

\end{document}